\documentclass[11pt]{article}
\usepackage{amsfonts,mathrsfs,amssymb,amsthm,mathptm}
\usepackage{amsmath,amscd}
\usepackage{mathptm,pslatex}
\usepackage{appendix}
\usepackage{multirow}
\usepackage{color}
\usepackage[all]{xy}
\usepackage{url}
\usepackage{cite}
\usepackage{exscale}
\usepackage{relsize}
\usepackage{bbm}
\begin{document}
\newtheorem{Def}{Definition}[section]
\newtheorem{Bsp}[Def]{Example}
\newtheorem{Prop}[Def]{Proposition}
\newtheorem{Theo}[Def]{Theorem}
\newtheorem{Lem}[Def]{Lemma}
\newtheorem{Koro}[Def]{Corollary}
\theoremstyle{definition}
\newtheorem{Rem}[Def]{Remark}

\newcommand{\add}{{\rm add}}
\newcommand{\con}{{\rm con}}
\newcommand{\gd}{{\rm gl.dim}}
\newcommand{\dm}{{\rm domdim}}
\newcommand{\tdim}{{\rm dim}}
\newcommand{\E}{{\rm E}}
\newcommand{\Mor}{{\rm Morph}}
\newcommand{\End}{{\rm End}}
\newcommand{\ind}{{\rm ind}}
\newcommand{\lcm}{{\rm lcm}}
\newcommand{\rsd}{{\rm res.dim}}
\newcommand{\rd} {{\rm rep.dim}}
\newcommand{\ol}{\overline}
\newcommand{\overpr}{$\hfill\square$}
\newcommand{\rad}{{\rm rad}}
\newcommand{\soc}{{\rm soc}}
\renewcommand{\top}{{\rm top}}
\newcommand{\stp}{{\mbox{\rm -stp}}}
\newcommand{\pd}{{\rm projdim}}
\newcommand{\id}{{\rm injdim}}
\newcommand{\fld}{{\rm flatdim}}
\newcommand{\fdd}{{\rm fdomdim}}
\newcommand{\Fac}{{\rm Fac}}
\newcommand{\Gen}{{\rm Gen}}
\newcommand{\fd} {{\rm findim}}
\newcommand{\Fd} {{\rm Findim}}
\newcommand{\Pf}[1]{{\mathscr P}^{<\infty}(#1)}
\newcommand{\DTr}{{\rm DTr}}
\newcommand{\cpx}[1]{#1^{\bullet}}
\newcommand{\D}[1]{{\mathscr D}(#1)}
\newcommand{\Dz}[1]{{\mathscr D}^+(#1)}
\newcommand{\Df}[1]{{\mathscr D}^-(#1)}
\newcommand{\Db}[1]{{\mathscr D}^b(#1)}
\newcommand{\C}[1]{{\mathscr C}(#1)}
\newcommand{\Cz}[1]{{\mathscr C}^+(#1)}
\newcommand{\Cf}[1]{{\mathscr C}^-(#1)}
\newcommand{\Cb}[1]{{\mathscr C}^b(#1)}
\newcommand{\Dc}[1]{{\mathscr D}^c(#1)}
\newcommand{\K}[1]{{\mathscr K}(#1)}
\newcommand{\Kz}[1]{{\mathscr K}^+(#1)}
\newcommand{\Kf}[1]{{\mathscr  K}^-(#1)}
\newcommand{\Kb}[1]{{\mathscr K}^b(#1)}

\newcommand{\modcat}{\ensuremath{\mbox{{\rm -mod}}}}
\newcommand{\Modcat}{\ensuremath{\mbox{{\rm -Mod}}}}
\newcommand{\stmodcat}[1]{#1\mbox{{\rm -{\underline{mod}}}}}
\newcommand{\pmodcat}[1]{#1\mbox{{\rm -proj}}}
\newcommand{\imodcat}[1]{#1\mbox{{\rm -inj}}}
\newcommand{\Pmodcat}[1]{#1\mbox{{\rm -Proj}}}
\newcommand{\Imodcat}[1]{#1\mbox{{\rm -Inj}}}
\newcommand{\PI}[1]{#1\mbox{{\rm -prinj}\,}}

\newcommand{\opp}{^{\rm op}}
\newcommand{\otimesL}{\otimes^{\rm\mathbb L}}
\newcommand{\rHom}{{\rm\mathbb R}{\rm Hom}\,}
\newcommand{\projdim}{\pd}
\newcommand{\Hom}{{\rm Hom}}
\newcommand{\Coker}{{\rm Coker}}
\newcommand{ \Ker  }{{\rm Ker}}
\newcommand{ \Cone }{{\rm Con}}
\newcommand{ \Img  }{{\rm Im}}
\newcommand{\Ext}{{\rm Ext}}
\newcommand{\StHom}{{\rm \underline{Hom}}}

\newcommand{\gm}{{\rm _{\Gamma_M}}}
\newcommand{\gmr}{{\rm _{\Gamma_M^R}}}

\def\vez{\varepsilon}\def\bz{\bigoplus}  \def\sz {\oplus}
\def\epa{\xrightarrow} \def\inja{\hookrightarrow}
\newcommand{\lra}{\longrightarrow}
\newcommand{\llra}{\longleftarrow}
\newcommand{\lraf}[1]{\stackrel{#1}{\lra}}
\newcommand{\llaf}[1]{\stackrel{#1}{\llra}}
\newcommand{\ra}{\rightarrow}
\newcommand{\dk}{{\rm dim_{_{k}}}}

\newcommand{\colim}{{\rm colim\, }}
\newcommand{\limt}{{\rm lim\, }}
\newcommand{\Add}{{\rm Add }}
\newcommand{\Tor}{{\rm Tor}}
\newcommand{\Cogen}{{\rm Cogen}}
\newcommand{\Tria}{{\rm Tria}}
\newcommand{\tria}{{\rm tria}}

{\Large \bf
\begin{center}
Centralizer matrix algebras and symmetric polynomials of partitions
\end{center}}
	
\medskip
\centerline{\textbf{Changchang Xi and Jinbi Zhang$^*$}}
	
\renewcommand{\thefootnote}{\alph{footnote}}
\setcounter{footnote}{-1} \footnote{$^*$Corresponding author's
		Email:  zhangjb@cnu.edu.cn}
\renewcommand{\thefootnote}{\alph{footnote}}
\setcounter{footnote}{-1} 
\renewcommand{\thefootnote}{\alph{footnote}}
\setcounter{footnote}{-1} 
	
\begin{abstract}
For a field $R$ of characteristic $p\ge 0$ and a matrix $c$ in the full $n\times n$ matrix algebra $M_n(R)$ over $R$, let $S_n(c,R)$ be the centralizer algebra of $c$ in $M_n(R)$. We show that $S_n(c,R)$ is a Frobenius-finite, $1$-Auslander-Gorenstein, and gendo-symmetric algebra, and that the extension $S_n(c,R)\subseteq M_n(R)$ is separable and Frobenius. Further, we study the isomorphism problem of invariant matrix algebras.
Let $\sigma$ be a permutation in the symmetric group $\Sigma_n$ and $c_{\sigma}$ the corresponding permutation matrix in $M_n(R)$. We give sufficient and necessary conditions for the invariant algebra $S_n(c_{\sigma},R)$ to be semisimple. If $R$ is an algebraically closed field, we establish a combinatoric characterization of when
two semisimple invariant $R$-algebras are isomorphic in terms of the cycle types of permutations.
\end{abstract}

{\footnotesize\tableofcontents\label{contents}}	
	
\section{Introduction}
Let $R$ be a unitary ring and $n$ a natural number. We set $[n]:=\{1, 2, \cdots, n\}$ and denote by $\Sigma_n$ the symmetric group of permutations on $[n]$, by $M_n(R)$ the ring of $n\times n$ matrices over $R$ and by $e_{ij}$ the matrix units of $M_n(R)$ for $i,j\in[n]$. The identity matrix in $M_n(R)$ is denoted by $I_n$.

For a nonempty set $C$ of matrices in $M_n(R)$, we define the \emph{centralizer algebra} of $C$ by
$$S_n(C,R):=\{a\in M_n(R)\mid ac=ca \,\mbox{ for all }\, c\in C\} $$
In case $C=\{c\}$, we write $S_n(c,R)$ for $S_n(\{c\},R)$.

Given a subgroup $G$ of $\Sigma_n$, we associate a set $C_G$ of all permutation matrices $c_{\sigma}$ with $\sigma\in G$, and define the $G$-\emph{invariant matrix ring} over $R$ of degree $n$ by $S_n(G,R):=S_n(C_G,R)$. We write $S_n(\sigma,R)$ for $S_n(c_{\sigma},R)$.

Invariant algebras and rings can be traced back to the classical invariant theory (see \cite{weyl}). If $C$ consists of nilpotent matrices and $R$ is an algebraically closed field, then the variety consisting of nilpotent matrices in $S_n(C,R)$ is of great interest in understanding properties of semisimple Lie algebras (see \cite{Pan, Pre}). Note that $S_n(\sigma, R)$ is a generalization of centrosymmetric matrix algebras studied in \cite{XY}. As is known, centrosymmetric matrices have significant applications in Markov processes \cite{W}, engineering problems and quantum physics \cite{Datta}. Centralizer algebras $S_n(c,R)$ include the Auslander algebra of $R[x]/(x^n)$ (see \cite{Xi2021}), which plays a crucial role in describing the orbits of parabolic subgroups acting on their unipotent radicals. Moreover, it is shown in \cite{Xi2021} that if $R$ is an algebraically closed field then $S_n(c,R)$ is a cellular algebra and the extension $S_n(c,R)\subseteq M_n(R)$ is a Frobenius extension. But over an arbitrary field, these results are still to be understood.
Recall that an extension $S\subseteq R$ of rings is called a \emph{Frobenius extension} if $_SR$ is a finitely generated projective $S$-module and $_RR_S\simeq \Hom_S(_SR,{}_SS)$ as $R$-$S$-bimodules. Such extensions have intriguing interest in many aspects in mathematics (see \cite{Kadison1999}).

For $\sigma\in\Sigma_n$,  we can write $\sigma$ as a product of disjoint cycles, say $\sigma=\sigma_1\sigma_2\cdots\sigma_s$, where $\sigma_i$ is a $\lambda_i$-cycle with $\lambda_i\ge 1$. Note that the product is unique up to re-ordering of these cycles. Thus we have a decomposition of $n=\lambda_1+\lambda_2+\cdots +\lambda_s$. By assuming $\lambda_1\ge \cdots\ge \lambda_s\ge 1$, we get a partition $\lambda=(\lambda_1,\cdots,\lambda_s)$ of $n$ with $s$ parts, which is called the \emph{cycle type} of $\sigma$. This numeric data seems to play a central role in understanding the structure of $S_n(\sigma,R)$.

This note has two purposes. The first one is to reveal some nice homological properties of $S_n(c,R)$ and to show that $S_n(c,R)\subseteq M_n(R)$ is always a Frobenius extension for arbitrary field $R$ and $c\in M_n(R)$, including $c$ being nilpotent. For the nilpotent case, one may see \cite[Theorem 1.1(2)]{cm} for some recent results. The second one is to give a combinatoric characterization in terms of cycle types of when two semisimple invariant matrix algebras are isomorphic. Here, the general question reads as follows.

\smallskip
{\it Let $R$ be a field  and $\sigma,\tau\in \Sigma_n$. What are the necessary and sufficient conditions for $S_n(\sigma,R)$ and $S_n(\tau,R)$ to be isomorphic in term of partitions?}

Recall that the dominant dimension of an $A$-module $_AM$ is the supremum of $n$ such that the first $n$ terms in a minimal injective resolution of $M$ are projective modules. An algebra $A$ over a field is called \emph{Gorenstein} if $\id(_AA)<\infty$ and $\id(A_A)$ $< \infty$; and $n$-\emph{Auslander-Gorenstein} if $\id(_AA)\le n+1\le \dm(_AA)$ (see \cite{Iyama2018}). Here, for an $A$-module $M$, we denote by $\id(_AM)$ and $\dm(_AM)$ for the injective and dominant dimensions of $_AM$, respectively.
It is noted that an algebra $A$ is $n$-Auslander-Gorenstein if and only if so is $A^{\opp}$ (see \cite[Proposition 4.1(b)]{Iyama2018}).
Thus Auslander-Gorenstein algebras are always Gorenstein.

Following \cite{hx}, we define the \emph{Frobenius part} of a finite-dimensional algebra $A$ over a field $k$ to be the endomorphism algebra of the projective module $Ae$ with $e^2=e$ such that $\add(Ae)$ is just the category of those projective $A$-modules that remain projective under any positive power of the Nakayama functor $\Hom_k(A,k)\otimes_A-$. Clearly, the Frobenius part of $A$ may not be basic, but unique up to Morita equivalence. The algebra $A$ is said to be \emph{Frobenius-finite} if its Frobenius part is representation-finite. Note that Frobenius parts of algebras play an important role in understanding both the Auslander-Reiten conjecture on stable equivalences and the lifting of stable equivalences to derived equivalences (see \cite{hx} for details).
Recall that a finite-dimensional algebra $A$ over a field $k$ is \emph{symmetric} if $A \simeq \Hom_k(A,k)$ as $A$-bimodules; and \emph{gendo-symmetric} \cite{FK} if it is isomorphic to the endomorphism algebra of a generator for a finite-dimensional symmetric $k$-algebra.

Our first main result reads as follows.

\begin{Theo}\label{thm1.1} If $R$ is a field and $c\in M_n(R)$, then

$(1)$ $S_n(c, R)\subseteq M_n(R)$ is a separable Frobenius extension.

$(2)$ $S_n(c, R)$ is a Frobenius-finite, $1$-Auslander-Gorenstein and gendo-symmetric algebra. Moreover, the Frobenius part of $S_n(c, R)$ is always a symmetric algebra.
\end{Theo}

To state our result on the isomorphism problem specifically, we introduce a few definitions. For a partition $\lambda=(\lambda_1,\lambda_2,\cdots,\lambda_s)$ of $n$ with $s$ parts, we define the so-called partition polynomial $\epsilon_{\lambda}(x)\in \mathbb{Z}[x]$, which is the elementary symmetric polynomials $g_i(\lambda)$ in $\lambda_1,\cdots,\lambda_s$, over the canonical bisemigroup of natural numbers
$$g_i(\lambda):=\sum_{1\le k_1<k_2<\cdots<k_i\le s}\lambda_{k_1}\centerdot\lambda_{k_2}\centerdot\cdots\centerdot\lambda_{k_i}, $$
where $\lambda_i\centerdot\lambda_j$ means the greatest common divisor of $\lambda_i$ and $\lambda_j$.
The \emph{partition polynomial} of $\lambda$ is defined as
$$ \epsilon_{\lambda}(x):= x^{s-1}- \frac{g_{s-1}(\lambda)}{g_{s}(\lambda)}x^{s-2}+\cdots + (-1)^{s-2}\frac{g_{2}(\lambda)}{g_s(\lambda)}x +(-1)^{s-1} \frac{g_{1}(\lambda)}{g_s(\lambda)}\in  \mathbb{Z}[x], $$
The coefficients of $\epsilon_{\lambda}(x)$ can be calculated graphically (see Section \ref{semirings} for details).

\begin{Theo}\label{iso-m} Let $R$ be a field of characteristic $p\ge 0$, let $\sigma\in \Sigma_n$ and $\tau \in \Sigma_m$ be of cycle types $\lambda=(\lambda_1,\lambda_2,\cdots,\lambda_s)$ and $\mu =(\mu_1,\mu_2,\cdots,\mu_t)$, respectively. Then the following hold.

$(1)$ $S_n(\sigma,R)$ is semisimple if and only if $p\nmid \lambda_i$ for all $1\le i\le s$.

$(2)$ Assume that the field $R$ is an algebraically closed such that $S_n(\sigma,R)$ and $S_m(\tau,R)$ are semisimple. Then $S_n(\sigma,R)\simeq S_m(\tau,R)$ if and only if $m=n$ and $\epsilon_{\lambda}(x)=\epsilon_{\mu}(x)$, where $\epsilon_{\lambda}(x)\in \mathbb{Z}[x]$ is the partition polynomial of $\lambda$. In particular,
$S_n(\sigma,R)$ and $S_m(\tau,R)$ are Morita equivalent if and only if $(-1)^sd_{\lambda}\, \epsilon_{\lambda}(1) = (-1)^td_{\mu}\, \epsilon_{\mu}(1)$, where $d_{\lambda}$ denotes the greatest common divisor of $\lambda_1,\lambda_2,\cdots, \lambda_s$.
\end{Theo}

Thus the isomorphism problem for semisimple, invariant matrix algebras over an algebraically closed field can be read off from the numerical values of the elementary symmetric polynomials evaluated at the corresponding partitions (see Section \ref{semirings} for more details).

This article is outlined as follows. Section \ref{sect2} is devoted to some basic facts on invariant algebras of matrices, while Section \ref{sect3} provides a proof of Theorem \ref{thm1.1}. Section 4 contributes to showing Theorem \ref{iso-m}. During the course of the proof, we introduce a so-called polynomial equivalence relation and triangular divisor matrices of partitions based on elementary symmetric polynomials of partitions. Also, graphical calculations of the coefficients of partition polynomials are illustrated.
At the end of this section, we formulate a few unsolved problems suggested by the main results.

\section{Centralizer matrix algebras\label{sect2}}

In this section, we recall basic notions and establish primary facts on invariant matrix algebras.
Throughout the paper, $R$ denotes a unitary ring and $Z(R)$ its center.

For a subgroup $G$ of $\Sigma_n$, the permutation matrices $c_{\sigma}:=e_{1,(1)\sigma}+e_{2,(2)\sigma}+\dots + e_{n,(n)\sigma}$ with $\sigma \in G$ satisfy $(c_{\sigma})_{ij}=\delta_{(i)\sigma,j}$ for $ i,j\in [n],$ $c_{\sigma}c_{\tau}=c_{\sigma\tau}$, and $c'_{\sigma}=c_{\sigma^{-1}}=c^{-1}_{\sigma}$, where $a'$ is the transpose of the matrix $a$ and $\delta_{ij}$ is the Kronecker symbol. Clearly, $G$ acts on $[n]$ and therefore on $M_n(R)$ by
$$ M_n(R)\times G \lra M_n(R), \; (a_{i j})\cdot g :=(a_{(i)g,(j)g}), \; (a_{i j})\in M_n(R), g\in G.$$
The set of all $G$-fixed points of this action in $M_n(R)$ is defined by
$$(M_n(R))^G :=\{a\in M_n(R)\mid a_{i j}=a_{(i)\sigma,(j)\sigma}, \; \forall \; \sigma\in G, i,j\in [n]\}=\bigcap_{\sigma\in G}S_n(\sigma, R).$$

\begin{Lem}\label{lem2.1} $(1)$ $S_n(G,R)= (M_n(R))^G$. In particular, $Z(R)[c_{\sigma}]\subseteq Z(S_n(\sigma,R))$, where $Z(R)[c_{\sigma}]$ consists of all polynomials in $c_{\sigma}$ with coefficients in $Z(R)$.

$(2)$ $S_n(G,R)$ is a subring of the matrix ring $M_n(R)$.

$(3)$ The transpose of matrices is an anti-automorphism of order $2$ of the ring $S_n(G,R)$.

$(4)$ If $G$ is the subgroup of $\Sigma_n$ generated by $\tau_1, \cdots, \tau_r\in \Sigma_n$, then $S_n(G,R)=\cap_{i=1}^r S_n(\tau_i,R).$
\end{Lem}

{\it Proof.} (1) This is already noticed in \cite{stuart}. For the convenience of the reader, we sketch here a proof which will be used in the sequel. By definition, $c_{\sigma}a=ac_{\sigma}$ if and only if $(c_{\sigma}a)_{ij}=(ac_{\sigma})_{ij}$ for all $i,j\in [n]$. But
$$\begin{array}{ll} 1) & (c_{\sigma}a)_{ij}=\sum_{p=1}^n(c_{\sigma})_{ip}a_{pj}=\sum_{p=1}^n\delta_{(i)\sigma,p}a_{pj}=a_{(i)\sigma,j},\\
2) & (ac_{\sigma})_{ij}=\sum_{p=1}^na_{ip}(c_{\sigma})_{pj}=\sum_{p=1}^n a_{ip}\delta_{(p)\sigma,j}=a_{i,(j)\sigma^{-1}},\\
3) & (c_{\sigma}ac_{\sigma^{-1}})_{ij}=a_{(i)\sigma,(j)\sigma}. \end{array}$$
Thus $c_{\sigma}a=ac_{\sigma}$ if and only if $a_{(i)\sigma,j}=a_{i,(j)\sigma^{-1}}$ for all $i,j$ if and only if $a_{ij}=a_{(i)\sigma,(j)\sigma}$ for all $i,j\in [n]$.

The other statements are clear from definition. $\square$

\smallskip
By Lemma \ref{lem2.1}(1), $S_n(\sigma,R) = S_n(\langle \sigma\rangle,R)$, where $\langle \sigma\rangle$ is the subgroup of $\Sigma_n$ generated by $\sigma$.
The next lemma shows that $S_n(\sigma,R)$ depends only upon the cycle type $\lambda$ of $\sigma$. So we may also write $S_n(\lambda, R)$ for $S_n(\sigma, R)$. The proof is similar to the one of \cite[Lemma 2.1]{Xi2021}, we leave it to the reader.

\medskip
\begin{Lem}\label{lem4.3} Let $\sigma, \tau\in \Sigma_n$, and let $d\in M_n(R)$ be an invertible matrix. Then

$(1)$ $S_n(\sigma^{-1}, R)=S_n(\sigma,R)$. In particular, $S_n(c_{\sigma},R)= S_n(c'_{\sigma},R)$.

$(2)$ $S_n(dc_{\sigma} d^{-1},R)\simeq S_n(\sigma,R)$ as rings.

$(3)$ $S_n(\sigma\tau,R)\simeq S_n(\tau\sigma,R)$ as rings.
\end{Lem}

\medskip
Let $\sigma\in \Sigma_n$ be of the cycle type $(\lambda_1, \cdots, \lambda_s)$, $G=\langle \sigma\rangle$, and $\Lambda:=\{(i,j)\mid 1\le i,j\le n\}$. Then $G$ acts on $\Lambda$ by $(i,j)\cdot \sigma=\big((i)\sigma,(j)\sigma\big)$. Let $O_{(i,j)}$ denote the $G$-orbit of $(i,j)$. Then, according to 3), a matrix $a\in S_n(\sigma, R)$ if and only if $a$ takes a constant value on each orbit $O_{(i,j)}.$ Here, we regard $a\in M_n(R)$ as a function from $\Lambda$ to $R$. Now we describe the number of $G$-orbits in $\Lambda.$

We assume $n\ge 2$ and consider the three cases.

(i) $\sigma=(12 \cdots n)$. In this case, $G$ has order $n$ and acts on $\Lambda$ by $(i,j)\sigma=(i+1,j+1)$. Thus the stabilizer subgroup of $G$ fixing $(i,j)\in \Lambda$ is trivial and the $G$-orbit of $(i,j)$ has $|G|$ elements. Hence there are $n$ $G$-orbits of $\Lambda$, they are
$\mathcal{O}_{(1,j)}$, $1\le j\le n.$

(ii) $\sigma=(12\cdots p)(p+1,\cdots, p+q)$ with $n=p+q$. In this case, $|G|=[p,q]$, the least common multiple of $p$ and $q$. Let $X:=\{(i,j)\mid 1\le i\le p, p+1\le j\le n\}$ with the $G$-action given by $(i,j)\sigma=(i+1,j+1).$ Let $\sigma_1=(1 2\cdots p)\in \Sigma_p$ and $\sigma_2=(p+1, \cdots, p+q)\in \Sigma_q$, here we consider $\Sigma_q$ as the symmetric group of permutations on $\{p+1, p+2,\cdots, p+q \}$. If $((i)\sigma^l, (j)\sigma^l)=(i,j)$ for a natural number $1\le l\le [p,q]$, then $((i)\sigma_1^l, (j)\sigma_2^l)=(i,j)$. This implies that $p|l$ and $q|l$. Thus $[p,q]|l$ and $l=[p,q]$. Hence the stabilizer subgroup of $G$ fixing $(i,j)$ is trivial and the $G$-orbit of $(i,j)$ has $[p,q]$ elements. Therefore the number of $G$-orbits of $X$ is $pq/[p,q]$.

(iii) $\sigma=\sigma_1\sigma_2\cdots\sigma_s$, a product of disjoint cycles with $s$ a natural number, where $\sigma_i$ is a cycle of length $\lambda_i\ge 1$ with the content $X_i$, the set of all numbers appearing in $\sigma_i$. Thus $|G|=[\lambda_1,\lambda_2,\cdots,\lambda_s]$ and $X_i$ is a $G$-orbit for $1\le i \le s$, where $[\lambda_1, \lambda_2,\cdots,\lambda_s]$ denotes the least common multiple of $\lambda_1,\lambda_2, \cdots, \lambda_s$. Therefore $X_i\cap X_j=\emptyset$ for $1\le i\ne j\le s,$ and $\cup_{j=1}^sX_j=\{1,2, \cdots, n\}.$ The set $\Lambda$ can be identified with the disjoint union
$$ \dot{\bigcup}_{1\le i,j\le s}X_i\times X_j.$$
For $(x,y)\in X_i\times X_j$, the $G$-orbit $\mathcal{O}_{(x,y)}$ of $(x,y)$ depends only upon $\sigma_i$ and $\sigma_j$. Thus, if $i=j$, then we are in case (i) and the number of $G$-orbits is $\lambda_i$. If $i\ne j$, then we are in the case (ii) and the number of $G$-orbits of $X_i\times X_j$ is $\lambda_i\lambda_j/[\lambda_i,\lambda_j]$. Thus
the number of $G$-orbits of $\Lambda$ is
$$ \ell_n(\sigma):=\sum_{1\le i,j\le s}\frac{\lambda_i\lambda_j}{[\lambda_i,\lambda_j]}=\sum_{1\le i,j\le s} \lambda_i\centerdot\lambda_j,$$ where $\lambda_i\centerdot\lambda_j$ stands for $\gcd(\lambda_i,\lambda_j)$, the greatest common divisor of $\lambda_i$ and $\lambda_j.$

We define $C(\lambda):=(\lambda_i\centerdot\lambda_j)_{s\times s}\in M_s(\mathbb{N}).$
This symmetric matrix is called the \emph{greatest common divisor matrix} of $\lambda$ in the literature (see \cite{BL,Li}). In this note, it will be called the \emph{dimension matrix} of $S_n(\sigma,R)$. Further, let
$$f_i:=\sum_{j\in X_i}e_{jj}, \; 1\le i\le s, \quad f_{ij}:=\sum_{(p,q)\in \mathcal{O}_{(i,j)}}e_{p,q}, \; (i,j)\in \Lambda.$$
Then $f_1, \cdots, f_s$ are pairwise orthogonal idempotent elements in $S_n(\sigma, R)$ and $\sum_{j=1}^sf_j=I_n$ in $M_n(R)$. This yields the matrix decomposition of $S_n(\sigma,R)$
$$ S_n(\sigma,R)=\begin{pmatrix}
{f}_1S_n(\sigma,R)f_1 & f_1S_n(\sigma,R)f_2 & \cdots & f_1S_n(\sigma,R)f_s\\
f_2S_n(\sigma,R)f_1& f_2 S_n(\sigma,R)f_2 & \cdots &f_2S_n(\sigma,R)f_s\\
\vdots & \vdots &\ddots & \vdots\\
f_{s}S_n(\sigma,R)f_1& f_s S_n(\sigma,R)f_2& \cdots & f_s S_n(\sigma,R)f_s
\end{pmatrix}.$$

The transpose of matrices is an $R$-involution on $S_n(\sigma,R)$ fixing $f_i$ for all $1\le i\le s$.

\begin{Lem}\label{lem21}
$(1)$ $S_n(\sigma, R)$ is a free $R$-module of rank $\ell_n(\sigma)$.

$(2)$ The $R$-rank of $f_iS_n(\sigma,R)f_j$ is $\lambda_i\centerdot\lambda_j.$

$(3)$ The dimension matrix $C(\lambda)$ of $S_n(\sigma,R)$ is positive definite if the numbers $\lambda_i$, $1\le i\le s$, are pairwise distinct. In this case, the determinant $\det C(\lambda)$ of $C(\lambda)$ satisfies $\varphi(\lambda_1)\cdots\varphi(\lambda_s)\le \det C(\lambda)\le \prod_{i=1}^s\lambda_i-\frac{s!}{2}$, where $\varphi(\lambda_i)$ is the Euler's totient function of $\lambda_i$.
\end{Lem}

{\it Proof.}
(1) and (2) are consequences of the definition of $S_n(\sigma,R)$ and the calculation of $G$-orbits of $\Lambda$ (see also \cite{sw}).
(3) is proved in \cite[Theorem 2]{BL} and \cite[Theorem 1 and Theorem 3]{Li}.
$\square$

\begin{Bsp}\label{bsp} $(1)$ {\rm $S_n((12\cdots n),R)\simeq R[C_n]$, where $R[C_n]$ is the group algebra of the cyclic group $C_n$ of order $n$ over $R$. Note that $c_{(12\cdots n)}=\sum_{j=1}^ne_{j,j+1}$ has order $n$ and $c_{(12\cdots n)}^n=f_1$ is the identity matrix. Here we understand $(n,n+1)=(n,1)$. The transpose of matrices in $S_n((12\cdots n), R)$ corresponds to the involution defined by $c_{(12\cdots n)}\mapsto c_{(12\cdots n)}^{-1}.$

$(2)$ If $\lambda=(\lambda_1,\cdots,\lambda_1)=(\lambda_1^s)$ where $\lambda_1$ appears $s$ times, then $S_{s\lambda_1}(\lambda,R)\simeq M_s(R[C_{\lambda_1}]).$
In fact, each element of $S_{s\lambda_1}(\lambda,R)$ can be partitioned as an $s\times s$ block-matrix such that each block has entries in $f_iS_{s\lambda_1}(\lambda,R)f_j\simeq R[C_{\lambda_1}]$. Thus $S_{s\lambda_1}(\lambda,R)\simeq M_s(R[C_{\lambda_1}]).$
}\end{Bsp}

By Lemma \ref{lem4.3}(2), up to isomorphism of rings, we may assume $$\sigma=(12\cdots\lambda_1)(\lambda_1+1,\cdots, \lambda_1+\lambda_2)\cdots (\lambda_1+\cdots+\lambda_{s-1}+1,\cdots, \lambda_1+\cdots +\lambda_s)$$ as a product of disjoint cycles. Let $\lambda=(\lambda_1,\cdots,\lambda_s)$, $X_1=\{1, \cdots, \lambda_1\}$, $X_i=\{j\mid \sum_{\ell=1}^{i-1}\lambda_{\ell} < j \le \sum_{\ell=1}^i\lambda_i\}$ for $2\le i\le s$, $d_{ij}:=\lambda_i\centerdot\lambda_j$.
Then
$$c_{\sigma}=\begin{pmatrix}
P_{\lambda_1} & 0            & \cdots & 0\\
  0           & P_{\lambda_2}& \cdots & 0\\
  \vdots       & \vdots       &\ddots  & \vdots\\
  0           & 0            & \cdots & P_{\lambda_s}
  \end{pmatrix}_{n\times n} \; \mbox{ where } \, P_m:=\begin{pmatrix}
0      &  1              &\cdots    & 0\\

\vdots &\vdots      & \ddots   & \vdots \\
  0    & 0           &   \cdots       &1 \\
  1    &  0             & \cdots         & 0\\
  \end{pmatrix}_{m\times m}$$

For $r_0, \cdots, r_{m-1}\in R$, we denote by $C(r_0,\cdots, r_{m-1})$ (or $C_m(r_0, \cdots,r_{m-1})$) the circulant matrix
$$\begin{pmatrix}
r_0          &  r_1        &r_2         &\cdots    & r_{m-1}\\
r_{m-1}      &  r_0        &r_1         &\cdots    & r_{m-2}\\
\vdots       &\ddots       & \ddots     & \ddots   & \vdots \\
 r_2         & \cdots      & \ddots     & \ddots   & r_1  \\
 r_1         & r_2         & \cdots     & r_{m-1}  & r_0\\
  \end{pmatrix}_{m\times m}=r_0I_m+r_1P_m+\cdots +r_{m-1}P_m^{m-1}.$$
Thus, for a cyclic group $C_m$ of order $m$, we have $$R[C_m]\simeq S_m((12\cdots m),R)=\{C(r_0,\cdots, r_{m-1})\mid r_0, \cdots, r_{m-1}\in R\}.$$

Recall that, given $\sigma\in \Sigma_n$ and $\tau\in \Sigma_m$, an $n\times m$ matrix $a=(a_{ij})\in M_{n\times m}(R)$ is called a $(\sigma,\tau)$-\emph{invariant matrix} if $c_{\sigma}a=ac_{\tau}$.

We see easily that an $s\times s$ block matrix $a=(A_{ij})$ with $A_{ij}\in M_{\lambda_i\times \lambda_j}(R)$ lies in $S_n(\sigma,R)$ if and only if $A_{ij}$ is a $(\sigma_i,\sigma_j)$-invariant matrix for each pair $1\le i,j\le s$. Thus
a matrix $a\in S_n(\sigma,R)$ can be written as a block-matrix form
$$ a=(A_{ij})=\begin{pmatrix}
A_{11} & A_{12} & \cdots & A_{1 s}  \\
A_{21} & A_{22} & \cdots & A_{2 s} \\
\vdots & \vdots & \ddots & \vdots \\
A_{s1} & A_{s2} & \cdots & A_{s s}
\end{pmatrix}_{s\times s} $$
where the diagonal $A_{ii}$ is a $\lambda_i\times \lambda_i$ circulant matrix and the off-diagonal $A_{ij}$ is a $\lambda_i\times \lambda_j$ matrix over $R$, having the block form:
$$ A_{ij}=\begin{pmatrix}
C(r_1,\cdots, r_{d_{ij}})  & \cdots & C(r_1,\cdots, r_{d_{ij}})  \\
\vdots & \cdots & \vdots  \\
C(r_1,\cdots, r_{d_{ij}})  & \cdots & C(r_1,\cdots, r_{d_{ij}})  \\
\end{pmatrix}_{\frac{\lambda_i}{d_{ij}} \times \frac{\lambda_j}{d_{ij}}} $$
for $r_1, \cdots, r_{d_{ij}}\in R$ (see \cite[p.260]{sw} for $R=\mathbb{C}$). In general, we have an isomorphism of $R$-modules
$$f_iS_n(\sigma,R)f_j\simeq \{A_{ij}\in M_{\lambda_i\times\lambda_j}(R)\mid r_1, \cdots, r_{d_{ij}}\in R\}.$$

\begin{Lem} \label{lem2.4} $(1)$ For  $1\le i\le s,$ $f_iS_n(\sigma,R)f_i\simeq R[C_{\lambda_i}]$ as rings.

$(2)$ If $\lambda_1=\lambda_2$, then $S_n(\sigma,R)f_1\simeq S_n(\sigma,R)f_2$ as left $S_n(\sigma,R)$-modules.

$(3)$ If $m:=n-\lambda_i$ and $f:=f_1+\cdots +f_{i-1}+ f_{i+1}+\cdots +f_s$, then $fS_n(\sigma,R)f\simeq  S_m(\sigma_1\cdots\sigma_{i-1}\sigma_{i+1}\cdots\sigma_s,R)$ as rings.
\end{Lem}

{\it Proof.} (1) If $\lambda_i=1$, then $f_iS_n(\sigma,R)f_i\simeq R$ and (1) holds. We assume $\lambda_i\ge 2$ and define $m_{i-1}:=\sum_{j=1}^{i-1}\lambda_j$, $g:=e_{m_{i-1}+1,m_{i-1}+2}+\cdots + e_{m_{i-1}+\lambda_i-1,m_{i-1}+\lambda_i}+e_{m_{i-1}+\lambda_i,m_{i-1}+1}$. Then
$g^{\lambda_i}=f_i$ and $g^j\ne f_i$ for all $1\le j<\lambda_i$. One can check that $\{f_i, g,g^2,\cdots, g^{\lambda_i-1}\}$ is an $R$-basis of $f_iS_n(\sigma,R)f_i$. Thus $f_iS_n(\sigma,R)f_i\simeq R[C_{\lambda_i}]$.

(2) By definition, $f_1=e_{11}+\cdots +e_{\lambda_1,\lambda_1}$ and $f_2=e_{\lambda_1+1,\lambda_1+1}+e_{\lambda_1+2,\lambda_1+2}+\cdots + e_{2\lambda_1,2\lambda_1}$. Let $a_0:=e_{1,\lambda_1+1}+e_{2,\lambda_1+2}+\cdots +e_{\lambda_1,2\lambda_1}$. Then $a_0\in S_n(\sigma,R)$ and $f_1a_0=a_0f_2=a_0.$ So, we define a map
$$ \psi: S_n(\sigma,R)f_1\lra S_n(\sigma,R)f_2, \quad af_1\mapsto af_1a_0 \mbox{ for } a\in S_n(\sigma,R).$$
Note that $\psi$ just moves the first $\lambda_1$ columns of $af_1$ entirely to the second $\lambda_1$ columns of $af_1a_0$, and sends other columns of $af_1$ to zero. Thus $\psi$ is a homomorphism of $S_n(\sigma,R)$-modules, and it is in fact an isomorphism of $S_n(\sigma,R)$-modules.

(3) We may assume $i=2$ since the argument below works for any $i$. If we delete the cycle $\sigma_2$ from $\sigma$ and consider $\tau:=\sigma_1\sigma_3\cdots\sigma_s$ as a permutation on
$Y:=[n]\setminus X_2=X_1\cup (X_3\cup\cdots\cup X_s)$, then $S_m(\sigma_1\sigma_3\cdots\sigma_s,R)$ is well defined.  Observe that $\sigma$ and $\tau$ have the same action on $Y$ and that a matrix $a$ in $fS_n(\sigma,R)f$ is of the form
$$\begin{pmatrix}
A & 0 & U  \\
0 & 0 & 0 \\
V & 0 & B \\
\end{pmatrix}_{n\times n} $$
where $A$ is a matrix indexed by $X_1\times X_1$, $U$ is a matrix indexed by $X_1\times (X_3\cup\cdots \cup X_s)$, $V$ is a matrix indexed by $(X_3\cup\cdots \cup X_s)\times X_1$ and $B$ is a matrix indexed by $(X_3\cup\cdots \cup X_s)\times (X_3\cup\cdots \cup X_s)$. Since $a\in fS_n(\sigma,R)f$, the block matrix
$$b:=\begin{pmatrix}
A  & U  \\
V  & B \\
\end{pmatrix}$$ is in $S_m(\tau,R)$. Conversely, given a matrix $b\in S_m(\tau,R)$, we partition $b$ into a block form according to $Y$, say the above form. Then we get an $n\times n$ matrix by putting $0$-blocks in the middle row and column, and this matrix is in $fS_n(\sigma,R)f$. Note that the identity in $fS_n(\sigma,R)f$ is $f=\begin{pmatrix}
1 & 0 & 0  \\
0 & 0 & 0 \\
0 & 0 & 1 \\
\end{pmatrix}.$ In this way, we get a one-to-one correspondence which preserves identity, addition and multiplication of matrices. This shows  $fS_n(\sigma,R)f\simeq S_m(\tau,R)$. $\square$

\section{ Proof of Theorem \ref{thm1.1}\label{sect3}}

This section is devoted to a proof of the statements of Theorem \ref{thm1.1}.
We start with the following proposition.

\begin{Prop}\label{thm-ext} Let $k\subseteq K$ be an extension of fields with $\dim_k(K)<\infty$, and let $B\subseteq A$ be an extension of finite-dimensional algebras over $k$. Then $B\otimes_kK\subseteq A\otimes_kK$ is a Frobenius extension if and only if $B\subseteq A$ is a Frobenius extension.
\end{Prop}

To prove Proposition \ref{thm-ext}, we begin with a few lemmas.

\begin{Lem} \label{ext-iso} $(1)$ Let $k\subseteq K$ be an extension of fields. Then it is  a Frobenius extension if and only if $\dim_k(K)$ is finite.

$(2)$ Let $A$ be an algebra over a field, and let $M$ and $N$ be finite-dimensional
$A$-modules. Then the following are equivalent for $M$ and $N$:

{\rm (i)} $M\simeq N$.

{\rm (ii)} There is an integer $n\ge 1$ such that $M^n\simeq N^n$.

\end{Lem}

{\it Proof.} (1) is trivial. (2) is not difficult, but for the convenience of the reader, we include here a proof. We show that (ii) implies (i). Suppose $M\simeq \bigoplus_{i=1}^sX_i^{s_i}$ and $N\simeq \bigoplus_{j=1}^tY_i^{t_j}$, where $X_i$ and $Y_j$ are indecomposable modules with $X_i\not\simeq X_p$, $Y_j\not\simeq Y_q$ for $1\le i\neq p\le s$ and $1\le j\ne q\le t$, $s_i\ge 1$ and $t_j\ge 1$. It follows from $M^n\simeq N^n$ that $s=t$ and, up to a permutation of these $Y_j$, that $X_i\simeq Y_i$  and $s_i=t_i$ for all $i$. Thus $M\simeq N$. $\square$

\medskip
For $k$-algebras $A$ and $B$, we may form the tensor product $A\otimes_kB$ of algebras over $k$. Given an $A$-module $_AX$ and a $B$-module $_BY$, the tensor product $X\otimes_kY$ has a left $(A\otimes_kB)$-module structure in a natural way, and the Hom-set $\Hom_k(X,Y)$ is a $B$-$A$-bimodule defined by $bf: x\mapsto b(xf)$, $fa: x\mapsto (ax)f$ for $x\in X, a\in A, b\in B$ and $f\in \Hom_k(X,Y)$.

The first statement (1) of the next lemma is known in \cite[Lemma 4]{sugano70} and \cite[Corollary 2.8]{hs68}, while the second statement (2) is in \cite[Chap. XI, Theorem 3.1, p. 209-210]{CE}. In fact, (1) can be proved by applying (2).

\begin{Lem} \label{CE-Sugano} $(1)$ Let $B\subseteq A$ be an extension of algebras over a field $k$ and $C$ be an algebra over $k$. If $B\subseteq A$ is a Frobenius (or separable) extension, then the extension $B\otimes_kC\subseteq A\otimes_kC$ of tensor products of algebras is a Frobenius (or separable) extension.

$(2)$ Let $A, B, C,D, E,F$ be algebras over a field $k$, and let $_AX,{}_AY$ be $A$-modules and $_BX',{}_BY'$ be $B$-modules. Then $\Hom_{A\otimes_kB}(X\otimes_kX', Y\otimes_kY')\simeq \Hom_A(X,X')\otimes_k\Hom_B(Y, Y')$ as $k$-spaces. Moreover, if $_AX_C$, $_AY_E$, $_BX'_D $ and ${}_BY'_F$ are bimodules, then the foregoing isomorphism is a $(C\otimes_kD)$-$(E\otimes_kF)$-bimodule homomorphism.
\end{Lem}

\begin{Lem}\label{matalg-iso} If $k\subseteq K$ is an extension of fields and $c\in M_n(k)$, then $S_n(c,K)\simeq S_n(c,k)\otimes_kK$ as $K$-algebras.
\end{Lem}

{\it Proof.} Observe that $a\lambda = \lambda a\in S_n(c,K)$ for $a\in S_n(c,k)$ and $\lambda\in K$. If we choose a basis of the $k$-space $K$, say $\{v_i\mid i\in I\}$, where $I$ is an index set (not necessarily finite), then any matrix $a=(a_{ij})\in M_n(K)$ can be expressed uniquely as $a= \sum_{p\in J} a^{(p)}v_p$, with $J$ a subset of $I$ and $a^{(p)}\in M_n(k)$ for $p\in J$. This is due to the fact that $\{v_iI_n\mid i\in I\}$ is a basis of the left $M_n(k)$-module $M_n(K)$, where $I_n$ is the identity matrix. Hence a matrix $a=(a_{ij})\in M_n(K)$ belongs to $S_n(c,K)$ if and only if $a^{(i)}\in S_n(c,k)$ for all $i\in J$. This yields that the restriction of the multiplication map $\mu: M_n(k)\otimes_kK\ra M_n(K)$ to $S_n(c,k)\otimes_kK$ gives rise to a surjective homomorphism $\mu': S_n(c,k)\otimes_kK\ra S_n(c,K)$ of $K$-algebras. Thus Lemma \ref{matalg-iso} follows from the commutative diagram
$$\xymatrix{
S_n(c,k)\otimes_kK \ar[r]^-{\mu'}\ar@{^{(}->}[d] & S_n(c,K)\ar@{^{(}->}[d]\\
M_n(k)\otimes_kK \ar[r]_-{\sim}^-{\mu} & M_n(K)\\
}$$ with $\mu$ clearly an isomorphism of $K$-algebras. $\square$

\medskip
{\bf Proof of Proposition \ref{thm-ext}.} It follows from Lemma \ref{CE-Sugano}(1) that $B\otimes_kK\subseteq A\otimes_kK$ is  a Frobenius extension  if $B\subseteq A$ is a Frobenius extension. Conversely, we postulate that $B\otimes_kK\subseteq A\otimes_kK$ is  a Frobenius extension. Then $_{B\otimes_kK}(A\otimes_kK)$ is a finitely generated projective $B\otimes_kK$-module and $$\Hom_{B\otimes_kK}(A\otimes_kK, _{}B\otimes_kK)\simeq {}_{A\otimes_kK}(A\otimes_kK)_{B\otimes_kK}.$$First, we show that $_BA$ is finitely generated. Suppose that $\{\lambda_i\mid i\in I\}$ is a $k$-basis of $k$-space $K$ with $\lambda_0=1$. Then every element $a$ in $A\otimes_kK$ can be uniquely expressed as $\sum_{i\in I_a}a_i\otimes\lambda_i$ with $a_i\in A$, where $I_a$ is a subset of $I$. Let $x_1, \cdots, x_m\in A\otimes_kK$ be generators of $A\otimes_kK$ over $B\otimes_kK$ and $x_i=\sum_{j\in I_i} a_{ij}\otimes_k\lambda_j$ with $a_{ij}\in A$, where $I_i$ is an index set. Pick an element $a\in A$ and consider $a\otimes 1$. We write
$a\otimes 1= \sum_{j=1}^m y_jx_j,$ where $y_j=\sum_{p\in J_j}b_p\otimes\lambda_p\in B\otimes_kK$, and $\lambda_p\lambda_j=\sum_{q\in I_{pj}}\mu_{p j q}\lambda_q$, where $\mu_{pjq}\in k$ and $I_{pj}$ is a finite subset of $I$. Then $a\otimes 1$ is a finite summation of elements of the form $\mu_{pjq}b_pa_{ij}\otimes \lambda_q$ for $q\in I_{pj}$. By the unique expression of elements in $A\otimes_kK$, we see that $a$ can be generated over $B$ by all $a_{ij}$ with $1\le i\le m$ and $j\in \cup_{t=1}^m I_t$. Thus $_BA$ is finitely generated. Since $A\otimes_kK$ is projective over $B\otimes_kK$ and $B\otimes_kK$ is projective over $B$, we see that $A\otimes_kK$ is projective over $B$. This implies that $_BA$ is projective over $B$. Now it follows from Lemma \ref{CE-Sugano}(2) that
$\Hom_{B\otimes_kK}(_BA_A\otimes_kK, {}_{B}B_B\otimes_kK)\simeq \Hom_B(_BA,{}_BB)\otimes_kK$. Thus we have an isomorphism $ _AA_B\otimes_kK\simeq {}_A\Hom_B(_BA_A, {}_BB)_B\otimes_kK$ as $(A\otimes_kK)$-$(B\otimes_kK)$-bimodules. Just considering the $A$-$B$-bimodule structure of this isomorphism, we get $(_AA_B)^{\dim_k(K)}\simeq \Hom_B(_BA_A, {}_BB)^{\dim_k(K)}$ as $A$-$B$-bimodules. By Lemma \ref{ext-iso}(2), we must have $_AA_B\simeq \Hom_B(_AA,{}_BB)$ as $A$-$B$-bimodules. Thus $B\subseteq A$ is a Frobenius extension. $\square$

\medskip
{\bf Proof of Theorem \ref{thm1.1}(1)}.  Suppose that $R$ is any field and $c\in M_n(R)$. Let $\bar{R}$ be an algebraic closure of $R$. So all eigenvalues of $c$ lie in $\bar{R}$. We consider the smallest subfield $K$ of $\bar{R}$ containing $R$ and all eigenvalues of $c$, that is, $K$ is obtained from $R$ by adding all eigenvalues of $c$. Thus $K$ is a finite extension of $R$, that is, $\dim_R(K)<\infty$. We may assume that $K$ is a splitting field of the characteristic polynomial of $c$. Thus $c$ is similar to a Jordan-block matrix by a matrix in $GL_n(K)$. It follows from \cite[Theorem 1.2(2)]{Xi2021} that $S_n(c,K)\subseteq M_n(K) $ is a separable Frobenius extension. By Proposition \ref{thm-ext} and Lemma \ref{matalg-iso}, we deduce that $S_n(c,R)\subseteq M_n(R)$ is a separable Frobenius extension. $\square$

Now  we turn to the \textbf{proof of Theorem \ref{thm1.1}(2)}.
Assume that $c$ is a Jordan-block matrix with the same eigenvalues.
More precisely, suppose
$$c={\rm diag}(J_{1}^{b_1},J_{2}^{b_2},\dots,J_{s}^{b_s})\in M_{n}(R),$$
where $J_i$, the Jordan block of size $\lambda_i$ with the eigenvalue $r\in R$, appears $b_i$ times for $1\leq i\leq s$ and $\lambda_1 \textgreater \lambda_2 \textgreater \cdots \textgreater \lambda_s$.

As in \cite[Section 2]{Xi2021}, we define
$$m_0:=0,\;m_i:=\sum_{p=1}^{i}b_p,\; n_{ij}:=j\lambda_i+\sum_{p=1}^{i-1}b_p\lambda_p,\;1\le i\le s,1\le j\le b_i.$$
Then $m_s$ is the number of Jordan blocks of $c$.
For each $i \in [m_s]$, let $g(i)$ be the smallest $g(i) \in [s]$ such that $i\le m_{g(i)}$,  and let  $h(i):=i-m_{g(i)-1}\in [b_{g(i)}]$ and $\theta_{ij}:=\min\{\lambda_{g(i)},\lambda_{g(j)}\} \mbox{ for } j\in [m_s].$

For each $i \in [m_s]$, we define
$$f_i:=\sum_{p=n_{g(i)h(i)}-\lambda_{g(i)}+1}^{n_{g(i)h(i)}}e_{pp},$$
that is, $f_i$ is the identity matrix corresponding to the $i$-th block in the identity matrix $I_n$.
Then $\{f_1,f_2,\cdots,f_{m_s}\}$ is a complete set of orthogonal primitive idempotents of $S_n(c,R)$ by \cite[Lemma 2.6(4)]{Xi2021}.

Let $\Lambda := S_n(c,R)$, $\Lambda_{i j}:=f_iS_n(c,R)f_j$ and $\tilde{\Lambda}_{ij} := \{a\in M_{\lambda_{g(i)}\times \lambda_{g(j)}}(R)\mid J_{g(i)}a=aJ_{g(j)}\}$.
Given $1\le i,j \le m_s$ and $1\le p\le \theta_{ij}$, we define $$F_{ij}^{p}:=	\sum_{v=1}^{p}e_{n_{g(i)h(i)}-\lambda_{g(i)}+p-v+1,n_{g(j)h(j)}-v+1.}$$
Then $\{F_{ij}^{p}\mid 1\le p\le \theta_{ij}\}$ is an $R$-basis of $f_i S_n(c,R) f_j$ and
$\{F_{ij}^{p}\mid 1\le i,j\le m_s,1\le p\le \theta_{ij}\}$ is an $R$-basis of $S_n(c,R)$ by \cite[Lemma 2.6(2)-(3)]{Xi2021}. Moreover, we have the following property.

\begin{Lem}\label{basis}
{\rm (\cite[Lemma 2.5]{Xi2021})}
If $1\le i,j,k,l\le m_s,1\le p\le \theta_{ik}$ and $1\le q \le \theta_{kj}$, then $$ F_{ik}^pF_{lj}^q
=\left\{
\begin{array}{ll}
0   & if \; k\neq l \mbox{ or } \;p+q-\lambda_{g(k)}\textless 1,\\ F_{ij}^{p+q-\lambda_{g(k)}}  & if\; k=l \mbox{ and } p+q-\lambda_{g(k)}\ge 1.
\end{array} \right.$$
\end{Lem}

\medskip
Given $1\le u\le m_s$, we have $f_u=F^{\lambda_{g(u)}}_{uu}\in \Lambda f_u$.
By Lemma \ref{basis}, $f_uF_{u1}^{\lambda_{g(u)}}=F^{\lambda_{g(u)}}_{u1}
=F^{\lambda_{g(u)}}_{u1}f_1\in f_u \Lambda f_1$.
So we can define a map
$$\alpha_u:\Lambda f_u\lra\Lambda f_1, \quad af_u\mapsto af_uF^{\lambda_{g(u)}}_{u1} \mbox{ for } a\in \Lambda.$$
Note that $\{F_{iu}^{p}\mid 1\le i\le m_s,1\le p\le \theta_{iu}\}$ is an $R$-basis of $\Lambda f_u$ and $F_{iu}^{p}F^{\lambda_{g(u)}}_{u1}=F^p_{i1}$ for $1\le i\le m_s,1\le p\le \theta_{iu}$ by Lemma \ref{basis}.
Then $\{F_{i1}^{p}\mid 1\le i\le m_s,1\le p\le \theta_{iu}\}$ is an $R$-basis of the image of $\alpha_u$.
Let $M_u$ be the image of $\alpha_u$.
Obviously, $\Lambda f_u$ and $M_u$ have the same dimensions as $R$-spaces.
Thus $\alpha_u:\Lambda f_u\to M_{u}$ is an isomorphism of $\Lambda$-modules. If $1\le u\le b_1$, then $\theta_{iu}=\theta_{i1}$ and $M_u=\Lambda f_u$. In this case, $\alpha_u:\Lambda f_u\to\Lambda f_1$ is an isomorphism of $\Lambda$-modules.

Similarly, it follows from $f_1F^{\lambda_{g(u)}}_{1u}=F^{\lambda_{g(u)}}_{1u}
=F^{\lambda_{g(u)}}_{1u}f_u$ that we may define a map
$$\beta_u:f_u \Lambda\lra f_1\Lambda, \quad f_ua\mapsto F^{\lambda_{g(u)}}_{1u}f_ua \, \mbox{  for } a\in \Lambda.$$
Let $N_u$ be the image of $\beta_u$.
Then $\{F_{1i}^{p}\mid 1\le i\le m_s,1\le p\le \theta_{ui}\}$ is an $R$-basis of $N_u$, and
$\beta_u: f_u\Lambda \to N_u$ is an isomorphism of right $\Lambda$-modules. In particular, if $1\le u\le b_1$ then $\theta_{ui}=\theta_{1i}$, $N_u= f_u \Lambda$ and $\beta_u$ is an isomorphism of right $\Lambda$-modules.

In order to calculate the injective dimension of $\Lambda$, we need the following lemma.
For convenience, let $\lambda_{s+1}=0$, $g(m_s+1)=s+1$, $M_{m_s+1}=0$ and $N_{m_s+1}=0$.  We denote by $D=\Hom_R(-,R)$ the usual duality.
\begin{Lem}\label{lempinj}
$\Lambda f_1/M_u\simeq D(f_1\Lambda/N_u)$ as $\Lambda$-modules for $1\le u\le m_s+1$.
In particular, $\Lambda f_1$ is projective-injective.
\end{Lem}

{\it Proof.}
For $1\le u\le b_1$, $M_u=\Lambda f_1$ and $N_u =f_1\Lambda$, and therefore $\Lambda f_1/M_u=D(f_1\Lambda/N_u)=0$.
Now suppose $b_1+1\le u\le m_s+1$.

We utilize the technique: If $_A X$ and $Y_A$ are finite-dimensional modules with a non-degenerate $R$-bilinear form $\langle -,-\rangle:X\times Y\ra R$ which is associative: $\langle ax,y\rangle=\langle x,ya\rangle$ for all $x\in X$, $y\in Y$, $a\in A$, then $X\simeq D(Y)$ as $A$-modules.

For simplicity, let $f:=f_1$ and
$$\rho:f\Lambda f\lra R[x]/(x^{\lambda_1}),\; \sum_{i=1}^{\lambda_1}r_iF^i_{11}\mapsto \sum_{i=1}^{\lambda_1}r_i\bar{x}^{\lambda_1-i} \; \mbox{ for } r_i\in R,$$ be the canonical isomorphism of algebras. We define two $R$-linear maps $\pi$ and $\epsilon$ as follows.
$$\pi:R[x]/(x^{\lambda_1})\lra R[x]/(x^{\lambda_1-\lambda_{g(u)}}), \sum_{i=0}^{\lambda_1-1}r_i\bar{x}^i\mapsto \sum_{i=0}^{\lambda_1-\lambda_{g(u)}-1}r_i\bar{x}^i \; \mbox{ for } r_i\in R,$$ $$\epsilon:R[x]/(x^{\lambda_1-\lambda_{g(u)}})\lra R,\; \sum_{i=0}^{\lambda_1-\lambda_{g(u)}-1}r_i\bar{x}^i\mapsto \sum_{i=0}^{\lambda_1-\lambda_{g(u)}-1}r_i \; \mbox{ for } r_i\in R,$$
where $\bar{x}$ denotes the coset of $x$ in the quotient rings. Note that $(\sum_{i=1}^{\lambda_1}r_iF^i_{11})\rho
=\sum_{i=1}^{\lambda_1}r_i\bar{x}^{\lambda_1-i}
=\sum_{j=0}^{\lambda_1-1}r_{\lambda_1-j}\bar{x}^j.$
Then $\eta:=\rho \pi \epsilon: f \Lambda f \ra R,\; \sum_{i=1}^{\lambda_1}r_iF^i_{11}
\mapsto
\sum_{i=0}^{\lambda_1-\lambda_{g(u)}-1}r_{\lambda_1-i}$, is an $R$-linear map. Further, we define
$$\langle -,-\rangle:\Lambda f/M_u\times f\Lambda/N_u\lra R,\langle bf+M_u,fa+N_u\rangle:=(fabf)\eta\; \mbox{ for } a, b\in\Lambda.$$
We point out that $\langle -,-\rangle$ is independent of the choice of representatives of cosets in $\Lambda f/M_u$ and $f\Lambda/N_u$.
Indeed, given $bf\in M_u$, we have $abf\in M_u$ for $a\in \Lambda$ since $M_u$ is a $\Lambda$-module.
Note that $\{F_{11}^{p}\mid 1\le p\le \lambda_{g(u)}\}$ is an $R$-basis of $f M_u$ and $(F_{11}^{p})\eta
=(\bar{x}^{\lambda_1-p})\pi\epsilon
=(0)\epsilon=0$ for $1\le p\le \lambda_{g(u)}$.
Thus $\langle bf+M_u, fa+N_u \rangle=(fabf)\eta=0$ for $a\in \Lambda$.
Similarly, for  $fa\in N_u$, we have $\langle bf+M_u, fa+N_u \rangle=(fabf)\eta=0$ for $b\in \Lambda$.
Thus $\langle -,-\rangle$ is well defined.

Clearly, $\langle -,-\rangle$ is an $R$-bilinear form and for $c\in\Lambda$, $$\langle cbf+M_u,fa+N_u\rangle=(facbf)\eta=\langle bf+M_u,fac+N_u\rangle.$$ Moreover, $\langle -,-\rangle$ is non-degenerate.

In fact, suppose $bf\in \Lambda f$ such that $\langle bf+M_u,fa+N_u\rangle=0$ for all $a\in \Lambda$. Then $(f\Lambda bf)\eta=0$. Since $(f\Lambda bf)\rho\pi$ is a left ideal in $R[x]/(x^{\lambda_1-\lambda_{g(u)}})$ and $\Ker(\epsilon)$ does not contain any nonzero left ideal of $R[x]/(x^{\lambda_1-\lambda_{g(u)}})$, we get $(f\Lambda bf)\rho\pi=0$. In the following, we prove $bf\in M_u$.

Since $\{F_{i1}^{p}\mid 1\le i\le m_s,1\le p\le \lambda_{g(i)}\}$ is an $R$-basis of $\Lambda f$, we write $bf=\sum_{i=1}^{m_s}\sum_{p=1}^{\lambda_{g(i)}}
r_{ip}F^p_{i1}\in \Lambda f$ with $r_{ip}\in R$.
For $1\le j\le u-1$, we have $F^{\lambda_{g(j)}}_{1j}\in f\Lambda$ and $F^{\lambda_{g(j)}}_{1j}bf
=\sum_{i=1}^{m_s}\sum_{p=1}^{\lambda_{g(i)}}
r_{ip}F^{\lambda_{g(j)}}_{1j}F^p_{i1}
=\sum_{p=1}^{\lambda_{g(j)}}
r_{jp}F^p_{11}\in f\Lambda bf$ by Lemma \ref{basis}.
It follows from $(f\Lambda bf)\rho\pi=0$ that
$$(F^{\lambda_{g(j)}}_{1j}bf)\rho \pi
=(\sum_{p=1}^{\lambda_{g(j)}}
r_{jp}F^p_{11})\rho \pi
=(\sum_{p=1}^{\lambda_{g(j)}}r_{jp}\bar{x}^{\lambda_1-p})\pi
=(\sum_{q=\lambda_1-\lambda_{g(j)}}^{\lambda_1-1}r_{j,\lambda_1-q}\bar{x}^q)\pi
=\sum_{q=\lambda_1-\lambda_{g(j)}}^{\lambda_1-\lambda_{g(u)}-1}r_{j,\lambda_1-q}\bar{x}^q
=0$$
in $R[x]/(x^{\lambda_1-\lambda_{g(u)}})$.
Thus $r_{j,\lambda_1-q}=0$ for $\lambda_1-\lambda_{g(j)}\le q\le \lambda_1-\lambda_{g(u)}-1$, that is, $r_{jp}=0$ for $\lambda_{g(u)}+1\le p\le \lambda_{g(j)}$.
Then $bf=\sum_{i=1}^{u-1}\sum_{p=1}^{\lambda_{g(u)}}
r_{ip}F^p_{i1}+\sum_{i=u}^{m_s}\sum_{p=1}^{\lambda_{g(i)}}
r_{ip}F^p_{i1}$.
Note that $$\{F_{i1}^{p}\mid 1\le i\le m_s,1\le p\le \theta_{iu}\}=\{F_{i1}^{p}\mid 1\le i\le u-1,1\le p\le \lambda_{g(u)}\}\cup \{F_{i1}^{p}\mid u\le i\le m_s,1\le p\le \lambda_{g(i)}\}$$ is an $R$-basis of $M_u$. Then $bf\in M_u$.

Similarly, if $fa\in f\Lambda$ such that $\langle bf+M_u,fa+N_u\rangle=0$ for all $b\in \Lambda$, then
$fa\in N_u$. Thus $\langle -,-\rangle$ is a non-degenerate $R$-bilinear form.

Now, we define another $R$-linear map $$\psi_u:\Lambda f/M_u\lra D(f\Lambda/N_u),\, bf+M_u\mapsto (fa+N_u\mapsto \langle bf+M_u,fa+N_u\rangle) \mbox{ for } a,b\in \Lambda.$$
Then  $(fa+N_u)(c(bf+M_u))\psi_u=\langle c(bf+M_u),fa+N_u\rangle=\langle bf+M_u,(fa+N_u)c\rangle
$ = $ ((fa+N_u)c)(bf+M_u)\psi_u $ = $(fa+N_u)(c(bf+M_u)\psi_u)$ for $a,b,c\in \Lambda$. This means that $\psi_u$ is a homomorphism of $\Lambda$-modules. By definition,  $\theta_{i u}=\theta_{u i}$ for $1\le i,j\le m_s$. It follows from \cite[Lemma 2.6(2)]{Xi2021} that
$$\dim_R(\Lambda f_i)=\sum_{1\le j\le
m_s}\dim_R(f_j\Lambda f_i)=\sum_{1\le j\le
m_s}\theta_{ij}=\sum_{1\le j\le
m_s}\theta_{ji}=\sum_{1\le j\le
m_s}\dim_R(f_i\Lambda f_j)=\dim_R(f_i \Lambda)$$ for
$1\leq i\leq m_s.$ Thus $\dim_R(\Lambda f_i)=\dim_R(f_i \Lambda)$ for $1\leq i\leq m_s$.
Due to $\dim_R(M_u)=\dim_R(\Lambda f_u)$ and $\dim_R(N_u)=\dim_R(f_u \Lambda)$, we then get
$\dim_R(M_u)=\dim_R(N_u)$ and $\dim_R(\Lambda f/M_u)=\dim_R(f\Lambda /N_u)$. Since $\psi_u$ is injective by the non-degenerative form, it is an isomorphism of $\Lambda$-modules.
$\square$

\begin{Lem}\label{inj} $\Lambda$ is $1$-Auslander-Gorenstein. In particular, $\Lambda$ is Gorenstein.
\end{Lem}

{\it Proof.} Suppose $s=1$.
Then $\Lambda\simeq M_{b_1}(R[x]/(x^{\lambda_1}))$ by \cite[Lemma 2.6(1)]{Xi2021}.
Thus $\Lambda$ is self-injective and $\id(_{\Lambda}\Lambda)=0<2<\infty=\dm(\Lambda)$.
Hence $\Lambda$ is $1$-minimal Auslander-Gorenstein.

Suppose $s\geq 2$.
By Lemma \ref{lempinj}, $\Lambda f_1\simeq D(f_1 \Lambda)$ is a projective-injective module. For $1\le u\le b_1$, we have $\Lambda f_u\simeq \Lambda f_1$ and $f_u\Lambda\simeq f_1\Lambda$. Then $\Lambda f_u$ is projective and injective for $1\le u\le b_1$.
For $b_1+1\leq u\leq m_s$, there are two exact
sequences of $\Lambda$-modules:
$$\qquad\qquad 0\lra \Lambda f_u\lraf{\alpha_u} \Lambda f_1\lraf{\alpha'_u}\Lambda f_1/M_u\lra 0,$$
$$ (*)\quad 0\lra f_u\Lambda \lraf{\beta_u} f_1\Lambda \lraf{\beta_u'} f_1\Lambda/N_u\lra 0.$$
By applying $D=\Hom_R(-,R)$ to the exact sequence $(*)$,
we get the exact sequence of $\Lambda$-modules :$$ 0\lra D(f_1\Lambda/N_u)\lraf{D(\beta_u')} D( f_1\Lambda)\lraf{D(\beta_u)} D(f_u\Lambda)\lra 0.$$
By Lemma \ref{lempinj}, we have $\psi_u:\Lambda f_1/M_u\lraf{\sim} D(f_1\Lambda/N_u)$ as $\Lambda$-modules. This gives rise to
an exact sequence
$$(\star)\quad 0\lra \Lambda f_u\lraf{\alpha_u}\Lambda f_1
\lraf{\alpha_u'\psi_uD(\beta_u')}
D(f_1\Lambda)\lraf{D(\beta_u)} D(f_u\Lambda)\lra 0$$
of $\Lambda$-modules. Again by Lemma \ref{lempinj}, $\Lambda f_1\simeq D(f_1 \Lambda)$ is projective-injective. Since $\Lambda f_1$ is indecomposable, $\alpha_u$ and $D(\beta_u')$ are injective envelopes of $\Lambda f_u$ and $D(f_1\Lambda/N_u)$, respectively. Therefore
the sequence $(\star)$ is a minimal injective resolution of $_{\Lambda}\Lambda f_u$.
Thus $\id(_{\Lambda}\Lambda f_u)=2$ for $b_1+1\le u\le m_s$ and $\id(_{\Lambda}\Lambda)=2$.

For $b_1+1\le u\le m_s$, $\Lambda f_u$ is not an injective $\Lambda$-module. Otherwise, it follows from the monomorphism $\alpha_u$ that $\Lambda f_u$ is a direct summand of $\Lambda f_1$, and this would mean $\Lambda f_u\simeq \Lambda f_1$, a contradiction to $\lambda_1 > \lambda_2 > \cdots > \lambda_s$.
Thus $\Lambda f_1$, up to isomorphism of $\Lambda$-modules, is the unique indecomposable, projective-injective $\Lambda$-module. By $(\star)$, $\dm(\Lambda f_u)=2$ for $b_1+1\le u\le m_s$ and $\dm(\Lambda)=2$.

Thus $\id(_{\Lambda}\Lambda)=2=\dm(\Lambda)$ and $\Lambda$ is $1$-Auslander-Gorenstein. $\square$

Since an indecomposable projective module which is stable under the Nakayama functor must be projective-injective, the proof of Lemma \ref{inj} shows that $\Lambda f_1$ is the only (up to isomorphism) indecomposable projective-injective $\Lambda$-module stable under arbitrary positive power of the Nakayama functor $D\Hom_{\Lambda}(-,\Lambda)$. Thus $\Lambda$ has the Frobenius part $f_1\Lambda f_1$ which is isomorphic to $R[x]/(x^{\lambda_1})$.

\begin{Lem}\label{sym}
Let $A$ be a finite-dimensional $k$-algebra over a field $k$, and let $k\subseteq K$ be an extension of fields.

$(1)$ If $A$ is a symmetric $k$-algebra, then $A\otimes_kK$ is a symmetric $K$-algebra.

$(2)$ Suppose that the extension $k\subseteq K$ is finite. If $A\otimes_kK$ is a symmetric $K$-algebra, then $A$ is a symmetric $k$-algebra.
\end{Lem}
{\it Proof.} Let $D_k:=\Hom_k(-,k)$ for a field $k$. Then we have  the isomorphisms of $A\otimes_kK$-bimodules.
\begin{align*}
(\diamond)\quad
D_K(A\otimes_kK)
&=\Hom_K(A\otimes_kK,K)\\
&\simeq \Hom_k(A,\Hom_K(K,K))
\quad (\mbox{by the adjunction isomorphism})\\
&\simeq \Hom_k(A,K)\simeq \Hom_k(A,k\otimes_kK)\\
&\simeq \Hom_k(A,k)\otimes_kK
\quad ( A \mbox{ is a finite-dimensional } k\mbox{-space})\\
&=D_k(A)\otimes_kK.
\end{align*}

(1) Suppose that $A$ is a symmetric $k$-algebra. Then there exists an isomorphism $\eta: {_A}A_A\ra {_A}D_k(A)_A$ of $A$-bimodules. Then the induced homomorphism $\eta\otimes 1: A\otimes_kK \ra D_k(A)\otimes_kK$ is an isomorphism of $A\otimes_kK$-bimodules. By $(\diamond)$, we have $A\otimes_kK \simeq D_K(A\otimes_kK )$ as $A\otimes_kK$-bimodules. Thus $A\otimes_kK$ is a symmetric $K$-algebra.

$(2)$ Suppose that the extension $k\subseteq K$ is finite and $A\otimes_kK$ is a symmetric $K$-algebra. Then $A\otimes_kK\simeq D_K(A\otimes_kK)$ as $A\otimes_kK$-bimodules. By $(\diamond)$, we get $A\otimes_kK\simeq D_k(A)\otimes_kK$ as $A\otimes_kK$-bimodules. Just considering the $A$-bimodule structure of this isomorphism, we then obtain $A^{\dim_k(K)}\simeq D_k(A)^{\dim_k(K)}$ as $A$-bimodules. By Lemma \ref{ext-iso}(2), we have $A\simeq D_k(A)$ as $A$-bimodules, that is, $A$ is a symmetric $k$-algebra.
$\square$

\medskip
Theorem \ref{thm1.1}(2) is a summary of the next two results.

\begin{Prop}\label{frob-part} Let $k$ be a field and $d\in M_n(k)$. Then
$S_n(d,k)$ is always $1$-Auslander-Gorenstein and Frobenius-finite. Moreover, the Frobenius part of $S_n(d,k)$ is a symmetric $k$-algebra.
\end{Prop}

{\it Proof.} Let $R$ be a splitting field of the characteristic polynomial of $d$ over $k$. Then all eigenvalues of $d$ lie in $R$ and $R$ is a finite extension of $k$. Thus $d$ is similar to a Jordan-block matrix $c$ in $M_n(R)$. We can assume $c=\mbox{ diag}(c_1, \cdots, c_t)\in M_n(R),$
where $c_i\in M_{n_i}(R)$ is a Jordan-block matrix with the same eigenvalue $r_i$ and $r_i\ne r_j$ for $1\le i\ne j\le t$.  By \cite[Lemma 2.1(1)]{Xi2021}, $S_n(d,R)\simeq S_n(c,R)$ as $R$-algebras.

Let $\Lambda_i:=S_{n_i}(c_i,R)$. Then there is an algebra isomorphism
$S_n(c,R)\simeq{\rm diag}(\Lambda_1,\Lambda_2,\dots,\Lambda_t)$ by \cite[Lemma 2.7(2)]{Xi2021}.
According to Lemma \ref{inj}, there holds $\id(_{\Lambda_i}\Lambda_i)\le 2\le \dm(\Lambda_i)$ for $1\le i\le s$.
It then follows from $\id(S_n(c,R))=\max\{\id(_{\Lambda_i}\Lambda_i)\mid 1\leq i\leq t\}$ and $\dm(S_n(c,R))=\min\{\dm(\Lambda_i)\mid 1\leq i\leq t\}$ that $\id(S_n(c,R))\le 2\le \dm(S_n(c,R))$. Hence $S_n(c,R)$ is $1$-Auslander-Gorenstein, and therefore $S_n(d,R)$ is $1$-Auslander-Gorenstein. By Lemma \ref{matalg-iso}, $S_n(d,R)\simeq S_n(d,k)\otimes_k R$ as rings.
It follows from \cite[Lemma 5]{mueller}
that $S_n(d,k)$ and $S_n(d,R)$ have the same self-injective and dominant dimensions. Thus $S_n(d,k)$ is $1$-Auslander-Gorenstein.

Let $D_k:=\Hom_k(-,k)$, $B:=S_n(d,k)$ and $\Lambda:= S_n(d,R) \simeq S_n(d,k) \otimes_kR = B\otimes_kR$. Suppose that $eBe$ is the Frobenius part of $B$ with $e^2=e\in B$, so that $Be$ is the direct sum of all non-isomorphic indecomposable projective-injective modules that remain projective under any positive power of the Nakayama functor $\nu_B:=D_k\Hom_B(-,B)$. We may assume $e\ne 0$ and prove that $\Lambda(e\otimes 1)$ is stable under the Nakayama functor $D_R\Hom_{\Lambda}(-,\Lambda)$. Indeed, since the $B$-module $Be$ is basic and the projective module $\nu_B(Be)$ remains projective under the functor $\nu_B^i$ for $i\ge 0$, we get $Be\simeq \nu_B(Be)$. Moreover, there are the following isomorphisms of $\Lambda$-modules:
\begin{align*}
D_R\Hom_{\Lambda}(\Lambda(e\otimes 1),\Lambda)
&\simeq D_R\Hom_{B\otimes_kR}(Be\otimes_kR,B\otimes_kR)\\
&\simeq D_R(\Hom_{B}(Be,B)\otimes_kR)
\quad \mbox{(by Lemma \ref{CE-Sugano} (2))}\\
&= \Hom_{R}(\Hom_{B}(Be,B)\otimes_kR,R)\\
&\simeq \Hom_{k}(\Hom_{B}(Be,B),R)\\
&\simeq \Hom_{k}(\Hom_{B}(Be,B),k\otimes_kR)\\
&\simeq \Hom_{k}(\Hom_{B}(Be,B),k)\otimes_kR \quad (\Hom_{B}(Be,B) \mbox{ is a finite-dimensional } k\mbox{-space})\\
&\simeq Be\otimes_kR\\
&\simeq \Lambda(e\otimes 1).
\end{align*}
Since $\Lambda$ is isomorphic to $\prod_{i=1}^t\Lambda_{i}$ as algebras, we may assume that $\epsilon_1,\epsilon_2,\cdots,\epsilon_t\in\Lambda$ are pairwise orthogonal central idempotents such that $1=\epsilon_1+\epsilon_2+\cdots+\epsilon_t$ and $\Lambda\epsilon_i\simeq \Lambda_i$ as algebras for $1\le i\le t$. In the following, we identify $\Lambda\epsilon_i$ and $\Lambda_i$. Then
$$\begin{array}{rl}
\epsilon_i\Lambda(e\otimes 1)& \simeq \epsilon_iD_R\Hom_{\Lambda}(\Lambda(e\otimes 1),\Lambda) \simeq D_R\Hom_{\Lambda}(\Lambda(e\otimes 1),\Lambda\epsilon_i) \\ & \simeq D_R\Hom_{\Lambda}(\epsilon_i\Lambda(e\otimes 1),\Lambda\epsilon_i) \simeq D_R\Hom_{\Lambda_i}(\epsilon_i\Lambda(e\otimes 1),\Lambda_i)\end{array}$$
as $\Lambda_i$-modules for $1\le i\le t$. On the other hand, by the proof of Lemma \ref{inj}, we may assume that $\Lambda_i f_{i1}$ is the indecomposable, projective-injective module which defines the Frobenius part of $\Lambda_i$ for $1\le i\le t$. Thus $\epsilon_i\Lambda(e\otimes 1)\in \add(\Lambda_i f_{i1})$. Moreover, for $1\le i\le t$, $\add(\epsilon_i\Lambda(e\otimes 1))=\add(\Lambda_i f_{i1})$ if and only if $\epsilon_i\Lambda(e\otimes 1)\neq 0$ if and only $\epsilon_i(e\otimes 1)\ne 0$. Thus $\epsilon_i\Lambda(e\otimes 1)$ defines the Frobenius part of $\Lambda_i$ if $\epsilon_i(e\otimes 1)\ne 0$. Let $I$ be the set of $1\le i\le t$ such that $\epsilon_i(e\otimes 1)\neq 0$. Then $\add(\Lambda(e\otimes 1))= \add( \bigoplus_{i=1}^t\epsilon_i\Lambda(e\otimes 1)) =\add(\bigoplus_{i\in I}\Lambda_{i} f_{i1})$. This shows that $\End_{\Lambda}(\Lambda(e\otimes 1))\simeq (e\otimes 1)\Lambda(e\otimes 1)= eBe\otimes_kR$ can also be regarded as the Frobenius part of the algebra $\prod_{i\in I}\Lambda_{i}$. Since the Frobenius part of $\prod_{i\in I}\Lambda_{i}$ is representation-finite and symmetric, the algebra $eBe\otimes_kR$ is representation-finite and symmetric. By \cite[Lemma 3.2]{jen82}, if $eBe\otimes_kR$ is representation-finite, then $eBe$ is representation-finite, namely $B$ is Frobenius-finite. By Lemma \ref{sym}, if $eBe\otimes_kR$ is a symmetric $R$-algebra, then $eBe$ is a symmetric $k$-algebra. $\square$

The following corollary is motivated by comments of an anonymous referee.

\begin{Koro}
$S_n(d,k)$ is a gendo-symmetric algebra for any field $k$ and $d\in M_n(k)$.
\end{Koro}
{\it Proof.} We keep the notation in the proof of Proposition 3.9. By Proposition 3.9, $\dm(B)\ge 2$. Thus we can pick a projective-injective faithful right $B$-module $P_B$ and consider $H:=\End_{B^{\opp}}(P)$. Then $_HP_B$ has a natural bimodule structure. It follows from $P\in\add(B_B)$ that $_H\End_{B^{op}}(P)\in \add({_H}\Hom_{B^{\opp}}(B,P))$, namely $_HP$ is a generator for the module category of the $k$-algebra $H$. We claim that the canonical homomorphism
$$\phi:B\ra \End_{H}(P),b\mapsto (p\mapsto pb),b\in B,p\in P$$
is an isomorphism of $k$-algebras. Indeed, it follows from $\dm(B)\ge 2$ that there is an injective resolution of $B_B$: $0\ra B_B\ra P^u\lra P^v$ with positive integers $u,v$. Since $P_B$ is injective, the induced sequence of $H$-modules:
$$\Hom_{B^{op}}(P^v,P)\lra \Hom_{B^{op}}(P^u,P)\lra \Hom_{B^{op}}(B,P) \lra 0$$ is exact. This implies that the following diagram is commutative and exact:
$$\xymatrix@C=2pc@R=2pc{
0\ar[r] &B_B\ar[r]\ar[d]^{\phi_B} &P^u\ar[d]^{\simeq}_{\phi_{P^u}}\ar[r] &P^v\ar[d]^{\simeq}_{\phi_{P^v}}\\
0\ar[r] &\Hom_{H}(\Hom_{B^{op}}(B,P),P)\ar[r] &\Hom_{H}(\Hom_{B^{op}}(P^u,P),P)\ar[r] &\Hom_{H}(\Hom_{B^{op}}(P^v,P),P)
}$$
where the vertical maps are defined canonically. Thus $\phi_B$ is an isomorphism of $B^{op}$-modules. Clearly $\Hom_{H}(\Hom_{B^{op}}(B,P),P)\simeq\Hom_H(_HP,{}_HP)=\End_H(P)$ and $\phi$ is an algebra isomorphism.

To prove that $B$ is a gendo-symmetric algebra, it remains only to show that $H$ is a symmetric $k$-algebra. Note that $P\otimes_kR$ is a projective-injective $\Lambda^{op}$-module and therefore $ (P\otimes_kR)\epsilon_i$ is a projective-injective $\Lambda_i^{op}$-module for $1\le i\le t$. By the proof of Lemma \ref{inj} for right modules, we see that $f_{i1}\Lambda_i $ is the only (up to isomorphism) indecomposable projective-injective $\Lambda_i^{op}$-module for $1\le i\le t$. Thus we have
$(P\otimes_kR)\epsilon_i\in \add(f_{i1}\Lambda_i )$. Moreover, for $1\le i\le t$, $\add( (P\otimes_kR)\epsilon_i)=\add(f_{i1}\Lambda_i)$ if and only if $(P\otimes_kR)\epsilon_i\neq 0$.
Let $J:=\{i\mid 1\le i\le t, (P\otimes_kR)\epsilon_i \neq 0\}$. Then $\add({_{\Lambda}}P\otimes_kR)= \add( \bigoplus_{i=1}^t (P\otimes_kR)\epsilon_i) =\add(\bigoplus_{i\in J}f_{i1}\Lambda_{i}) $ and $\End_{\Lambda^{op}}(P\otimes_kR)$ is Morita equivalent to $\End_{\Lambda^{op}}(\bigoplus_{i\in J} f_{i1}\Lambda_{i})$, but the latter is isomorphic to $\prod_{i\in J}f_{i1}\Lambda_{i} f_{i1}$ as algebras. By \cite[Lemma 2.6(1)]{Xi2021}, $f_{i1}\Lambda_i f_{i1}$ is a symmetric $R$-algebra for $i\in J$, and therefore $\prod_{i\in J}f_{i1}\Lambda_{i} f_{i1}$ is a symmetric $R$-algebra. It is known that a finite-dimensional algebra over a field Morita (or derived) equivalent to a symmetric algebra is itself symmetric. Thus $\End_{\Lambda^{op}}(P\otimes_kR)$ is a symmetric $R$-algebra. Now, it follows from the $R$-algebra isomorphisms  $$H\otimes_k R=\End_{B^{op}}(P)\otimes_k R\simeq \End_{(B\otimes_kR)^{op}}(P\otimes_k R)=\End_{\Lambda^{op}}(P\otimes_kR)$$  that $H\otimes_k R$ is a symmetric $R$-algebra. Since the extension $k\subseteq R$ is finite, $H$ is a symmetric $k$-algebra by Lemma \ref{sym}. Hence $B$ is a gendo-symmetric $k$-algebra.
$\square$
	
\section{Proof of Theorem \ref{iso-m}}
In this section we introduce an equivalence relation on the set of partitions of $n$ with exactly $s$ parts by employing elementary symmetric polynomials. We then study combinatorics of this relation by a matrix norm over the semiring of nonnegative integers. Finally, we prove Theorem \ref{iso-m}.

\subsection{A new equivalence relation on partitions \label{semirings}}

Recall that a nonempty set $A$ together with two associative binary operations $+$ and $\cdot$, named addition and multiplication, respectively, is called a \emph{bisemigroup}. A bisemigroup $(A,+,\cdot)$ is called a \emph{semiring} if $(A,+)$ is commutative and the distributive laws hold: $x \cdot ( y + z) = x\cdot y + x\cdot z$ and $(x+y)\cdot z = x\cdot z + y\cdot z$ for $x, y, z\in A$. A \emph{commutative} semiring is a semiring such that the multiplication is commutative.

\begin{Def}  A bisemigroup $(A,+,\centerdot)$ is called a quarter-ring  if

$(1)$ $(A,+)$ and $(A,\centerdot)$ are commutative,

$(2)$ $a\centerdot a= a$ for $a\in A$, and

$(3)$ $(a+b)\centerdot a=b\centerdot a$ for $a,b\in A$.
\end{Def}

In a bisemigroup $(A, +, \cdot)$, for any positive integer $m$ and $a_i, b, d\in A$, we write $\sum_{i=1}^m a_i$ for $a_1+\cdots +a_m$, particularly, $ma$ for $\sum_{i=1}^m a$, and $a^{\cdot m}$ for $\underbrace{a\cdot a\cdot \cdots \cdot a}_m$, the product of $m$ copies of $a$.

In a quarter-ring $(A,+,\centerdot)$, the following hold.

$\bullet\quad$ If $a\centerdot b=d$, then $a\centerdot d=d.$ In particular, if $ a_i\centerdot b=d$ for $1\le i\le n$, then $a_1\centerdot\cdots\centerdot a_n\centerdot b=d.$

$\bullet\quad$ $(ma+b)\centerdot a=a\centerdot b $ and $(ma)\centerdot a= a$ for any positive integer $m$.

$\bullet\quad$ If $(A,+)$ is additionally a monoid with zero element $0$, then $0\centerdot a = a$ for $a\in A$.

$\bullet\quad$ $ a_1\centerdot a_2\centerdot\cdots\centerdot a_n=(a_1\centerdot a_2)\centerdot(a_2\centerdot a_3)\centerdot\cdots \centerdot (a_{n-1}\centerdot a_n).$

$\bullet\quad$ $ (a_1\centerdot a_2\centerdot\cdots\centerdot a_n)\centerdot b =(a_1\centerdot b)\centerdot (a_2\centerdot \cdots\centerdot a_n) =(a_1\centerdot b)\centerdot (a_2\centerdot b)\centerdot (a_3\centerdot \cdots\centerdot a_n)=\cdots =(a_1\centerdot b)\centerdot (a_2\centerdot b)\centerdot \cdots\centerdot (a_n\centerdot b).$

\begin{Bsp}\label{ex1.2}{\rm
(1) Let $\mathbb{N}$ be the set of natural numbers including $0$. Then $\mathbb{N}$ with the usual addition and multiplication forms a semiring.
The polynomial semiring $\mathbb{N}[x_1,\cdots,x_s]$ over $\mathbb{N}$ in variables $x_1, \cdots, x_s$ is the set of all polynomials over $\mathbb{N}$ with the usual addition and multiplication of polynomials. Similarly, we define the matrix semiring $M_n(R)$ over a semiring $R$.

(2) Let $a\centerdot b:= \gcd(a,b)$ for $a,b\in \mathbb{N}$. Then $(\mathbb{N},+,\centerdot)$ is a quarter-ring, but not a semiring.  It is called the \emph{canonical quarter-ring} of natural numbers.
}\end{Bsp}

Recall that a \emph{multiset} is a collection of elements possibly containing duplicates. For a multiset $\lambda=\{\lambda_1,\cdots,\lambda_s\}$ of $s$ elements in $\mathbb{N}$, we define a map $\varphi_{\lambda}$
from the polynomial semiring $\mathbb{N}[x_1, \cdots, x_s]$ to $(\mathbb{N},+,\centerdot)$ as follows:
$$\varphi_{\lambda}: \mathbb{N}[x_1,\cdots,x_s]\lra \mathbb{N}, \qquad \sum_{i_1,\cdots, i_s}a_{i_1,\cdots,i_s}x_1^{i_1}\cdots x_s^{i_s}\mapsto \sum_{i_1,\cdots, i_s}a_{i_1,\cdots,i_s}\lambda_1^{i_1}\centerdot\cdots\centerdot\lambda_s^{i_s}.$$
Note that $\varphi_{\lambda}$ is additive and can be extended to the matrix semiring over $\mathbb{N}[x_1,\cdots,x_s]$:
$$ \|-\|_{_{\lambda}}: M_n(\mathbb{N}[x_1,\cdots,x_s])\lra \mathbb{N},\qquad (a_{ij})_{n\times n}\mapsto \|(a_{ij})\|_{_{\lambda}}:=\sum_{1\le i,j\le n}(a_{ij})\varphi_{\lambda},$$
where $(a_{ij})\varphi_{\lambda}$ denotes the image of $a_{ij}$ under $\varphi_{\lambda}.$ The map $\|-\|_{_{\lambda}}$ is called the \emph{norm map}.

\medskip
Now we recall the following fact in combinatorics \cite[Theorem 6.1.1, Corollary 6.1.2, p. 614-615]{Brualdi2010}.

\begin{Lem}\label{in-ex}
Let $V_i$ be a finite set for $1\le i\le s$.
We define
$$V:=\bigcup_{1\le j \le s}V_j,\;\; g_i:=
\sum_{1\le k_1<k_2<\cdots<k_i\le s}\big|V_{k_1}\cap \cdots \cap V_{k_i}\big|,\; 1\le i\le s,$$
where $|S|$ indicates the number of elements of a set $S$.
Then

$(1)$ For $1\le i\le s$, the number $h_i$ of the elements of $V$ which are members of exactly $i$  sets of $V_1,V_2,\cdots,V_s$ is given by
$h_i=\sum_{k=0}^{s-i}(-1)^k C_{i+k}^i g_{i+k}$, where $C_{i+k}^i := \binom{i+k}{i}=\frac{(i+k)!}{i!k!}$ is the number of $i$-subsets of an $(i+k)$-element set.

$(2)$ Principle of Inclusion-Exclusion:
$\big|V\big|=\sum_{i=1}^{s}h_i=\sum_{i=1}^{s}(-1)^{i+1}g_i.$ Moreover, $\sum_{i=1}^s |V_i|= \sum_{i=1}^{s}ih_i.$
\end{Lem}

\smallskip
Let $R$ be an algebraically closed field of characteristic $p\ge 0$, and let $\sigma\in \Sigma_n$ be of cycle type $(\lambda_1,\lambda_2,\cdots,\lambda_s)$.
Suppose $p\nmid \lambda_i$ for $1\le i\le s$.
For $1\le i\le s$, we denote by $V_i$ the set of all $\lambda_i$-th roots of unity in $R$ and define $V:=\bigcup_{1\le j \le s}V_j.$
It follows from $p\nmid \lambda_i$ that $V_i$ is a cyclic group of order $\lambda_i$. Moreover, $$V_{k_1}\cap \cdots \cap V_{k_i}=\{\omega\in R\mid \omega^d=1\}$$ for $1\le k_1<k_2<\cdots<k_i\le s$, where $d:=\gcd(\lambda_{k_1},\lambda_{k_2},\cdots,\lambda_{k_i})$ is the greatest common divisor of $\lambda_{k_j},1\le j\le i$. Thus $\big|V_{k_1}\cap \cdots \cap V_{k_i}\big|=\gcd(\lambda_{k_1},\lambda_{k_2},\cdots,\lambda_{k_i})$ for $1\le k_1<k_2<\cdots<k_i\le s.$
We further define
$$g_i(\sigma):=\sum_{1\le k_1<k_2<\cdots<k_i\le s}\gcd(\lambda_{k_1},\lambda_{k_2},\cdots,\lambda_{k_i}),$$
and $H_i(\sigma)$ to be the subset of $V$ consisting of all elements which belong to exactly $i$ sets of the given $V_1,V_2,\cdots,V_s$, namely $$H_i(\sigma)=\{v\in V\mid \exists\, 1\le k_1< k_2<\cdots <k_i\le s, v\in \bigcap_{1\le j\le i}V_{k_j}, \, v\not\in \bigcup_{r\ne k_j,1\le j\le i}V_r\}.$$
Clearly, $H_i(\sigma)\cap H_j(\sigma)=\emptyset$ if $1\le i\ne j\le s$. By Lemma \ref{in-ex}(1), $h_i(\sigma):=|H_i(\sigma)|=\sum_{k=0}^{s-i}(-1)^k C_{i+k}^i g_{i+k}(\sigma).$
Obviously, $g_i(\sigma)$ and $H_i(\sigma)$ depend merely on the cycle type of $\sigma$. So $g_i(\lambda)$ and $h_i(\lambda)$ are well-defined functions from the set of all partitions of positive integers to $\mathbb{N}$ if we define $g_i=0$ and $h_i=0$ for $i> s$.

As in Example \ref{ex1.2}, we write $a\centerdot b$ for $\gcd(a,b)$ for $a,b\in \mathbb{Z}$. Then $g_i(x_1,\cdots, x_s)$ is just the $i$-th elementary symmetric polynomial over the canonical quarter-ring of natural numbers:
$$ g_i(x_1,\cdots, x_s)=\sum_{1\le j_1<j_2< \cdots <j_i\le s}x_{j_1}\centerdot x_{j_2}\centerdot\cdots \centerdot x_{j_i}$$
for $1\le i\le s$. We set $g_0(x_1,\cdots,x_s)=g_i(x_,\cdots,x_s)=0$ for $i > s$. Clearly, $g_i(x_1,\cdots,x_s)$ does not depend on any order of $x_1\ge x_2\ge \cdots\ge x_s$, and $g_s(x_1,\cdots,x_s)\mid g_i(x_1,\cdots,x_s)$ for $1\le i\le s$ and all $x_i\in \mathbb{N}$.

Recall that a \emph{partition} $\lambda$ of a positive integer $n$ is an $s$-tuple $(\lambda_1,\cdots, \lambda_s)$ of integers satisfying $\lambda_1\ge \lambda_2\ge \cdots\ge \lambda_s\ge 1$ and $\sum_{i=1}^s\lambda_i=n$. We write $\lambda=(\lambda_1,\cdots,\lambda_s)\vdash n$ for simplicity. Each $\lambda_i$ is called a \emph{part} of $\lambda$. Sometimes it is also convenient to think of partitions as multisets and to write $\lambda$ in the form $(\lambda_1^{a_1}, \lambda_2^{a_2},\cdots, \lambda_t^{a_t})$ with $\lambda_1>\lambda_2>\cdots>\lambda_t\ge 1$, where $\lambda_i$ appears $a_i\ge 1$ times in the multiset of $\lambda$.

Let $P(n)$ be the set of all partitions of $n$, $P(s,n)$ the set of partitions of $n$ with exactly $s$ parts, and $P_s^*(n)$ the set of partitions of $n$ with the largest part $s$. Then $P(s,n)$ and $P_s^*(n)$ have the same cardinality. This can be seen from Ferrers graphs (or Young diagrams) of partitions.

Convention: Given partitions $\lambda\in P(n)$ and $\mu\in P(m)$, we write $(\lambda,\mu)$ for the partition $\gamma\in P(n+m)$ such that the parts of $\gamma$ is the disjoint union of the parts of $\lambda$ with the parts of $\mu$. For example, if $\lambda=(4,3,2,1)$ and $\mu=(5,4,4,2,2,1,1)$, then $(\lambda,\mu)=(5,4,4,4,3,2,2,2,1,1,1)$. For $1\le d\in \mathbb{N}$, we define $d\lambda=(d\lambda_1,\cdots,d\lambda_s)\in P(s, dn)$ and $\lambda_{-} :=(\lambda_1,\cdots,\lambda_{s-1}) =(\lambda_1,\cdots,\lambda_{s-1},\hat{\lambda}_s)\in P(s-1,n-\lambda_s)$. Inductively, $\lambda_{-j}=(\lambda_1, \cdots, \lambda_{s-j})$. If $j\ge s$, we define $\lambda_{s-j}=(0,\cdots, 0)$. Then
$g_i(d\lambda) = d\, g_i(\lambda) $, and $$g_i(\lambda) = g_i(\lambda_-)\; +\sum_{1\le j_1<\cdots<j_{i-1} \le s-1}\lambda_{j_1}\centerdot\cdots \centerdot\lambda_{j_{i-1}}\centerdot\lambda_s \le g_i(\lambda_-)\; +\sum_{1\le j_1<\cdots<j_{i-1} \le s-1}\lambda_{j_1}\centerdot\cdots \centerdot\lambda_{j_{i-1}},$$ and therefore, if $\lambda_s=1$, then $g_i(\lambda)=g_i(\lambda_{-})+C_{s-1}^{i-1}.$ In general, we have
$$  g_i(\lambda_-) + C^{i-1}_{s-1} \le g_i(\lambda)
\le g_i(\lambda_-) + \min\{g_{i-1}(\lambda_-), \lambda_sC_{s-1}^{i-1}\}.$$

Now, we introduce a new equivalence relation $\sim$ on $P(s,n)$ and define a polynomial $\epsilon_{\lambda}(x)$ for each $\lambda\vdash n$.

\begin{Def}\label{character} $(1)$ The partition polynomial $\epsilon_{\lambda}(x)$ of $\lambda\in P(s,n)$ is
$$\epsilon_{\lambda}(x) :=  \sum_{i=0}^{s-1}(-1)^{s-1-i}\frac{g_{i+1}(\lambda)}{g_s(\lambda)}x^{i} = x^{s-1}- \frac{g_{s-1}(\lambda)}{g_{s}(\lambda)}x^{s-2}+\cdots + (-1)^{s-2}\frac{g_{2}(\lambda)}{g_s(\lambda)}x +(-1)^{s-1} \frac{g_{1}(\lambda)}{g_s(\lambda)}\in \mathbb{Z}[x]. $$

$(2)$ The equivalence relation $\sim$ on $P(s,n)$ (or $P(n))$ is given by
 $$ \lambda \sim \mu \; \mbox{ if  } \; \epsilon_{\lambda}(x)=\epsilon_{\mu}(x).$$
Equivalently, $\lambda\sim \mu$ if and only if $g_i(\lambda)g_s(\mu)=g_s(\lambda)g_i(\mu) \mbox{ for } 1\le i\le s$. This equivalence relation is called the \emph{polynomial equivalence} of partitions in $P(n)$.
 \end{Def}

\begin{Bsp}{\rm  (1) $\epsilon_{(1^n)}(x)=\sum_{i=0}^{n-1}(-1)^{n-1-i}C^{i+1}_nx^i,$ while $\epsilon_{(m,1^{n-m})}(x)=\sum_{i=1}^{n-m}(-1)^{n-m-i}C_{n-m+1}^{i+1}x^i +(-1)^{n-m}n$ for $1\le m<n$ and $\epsilon_{(n)}(x)=1$.

(2) In $P(2,p)$ with $p$ a prime, any two-part partitions are polynomial equivalent because the partition polynomial of each such partition is of the form $x-p$. In $P(3,11)$, there hold $(8,2,1)\sim (7,2,2)\sim (6,4,1)\sim (5,4,2)$, but $(8,2,1)\nsim (6,3,2)$.
}\end{Bsp}

Remark that $(-1)^{s-1}g_s(\lambda)\epsilon_{\lambda}(1)$ counts the number of distinct eigenvalues of $c_{\sigma}$ in $\mathbb{C}$, where $\sigma\in \Sigma_n$ has the cycle type $\lambda$ (see \cite[Lemma 6.3]{sw}).
Next, we consider when two partitions are polynomial equivalent. We first give a matrix interpretation of the polynomial equivalence.

\begin{Def} Given $\lambda=(\lambda_1,\cdots,\lambda_s)\in P(s,n)$, we associate it with a matrix $\Lambda_{\lambda}\in M_s(\mathbb{N})$ by setting
$$ \Lambda_{\lambda}:=\begin{pmatrix}
0& d_{12} & d_{13} &\dots & d_{1s}\\
0&                              0 & d_{23}   &\dots & d_{2s}\\		
\vdots&    \vdots                 &\ddots                          &\ddots&\vdots\\
0     &    0                      & \cdots                     &\ddots &  d_{s-1,s}\\
0	  &                        0  & 0                              &\dots &0\\
\end{pmatrix}$$ where $d_{ij}=\lambda_i\centerdot \lambda_j$ is the product in the canonical quarter-ring of natural numbers. This matrix is called the \emph{triangular divisor matrix} of $\lambda$. \end{Def}

Now we define a norm $\|\Lambda_{\lambda}\|$ of $\Lambda_{\lambda}$. First, we consider $\lambda_j, 1\le j\le s,$ as variables and then regard $\Lambda_{\lambda}$ as a matrix with the entries in the polynomial semiring $\mathbb{N}[\lambda_1,\cdots,\lambda_s]$, thinking of $\lambda_i\centerdot\lambda_j$ as the usual monomial $\lambda_i\lambda_j$ in $\mathbb{N}[\lambda_1,\cdots,\lambda_s]$. By applying $\|-\|_{_{\lambda}}$ to this matrix, we get the norm $\|\Lambda_{\lambda}\|$ of the triangular divisor matrix $\Lambda_{\lambda}$. That is, $\|\Lambda_{\lambda}\|=\sum_{1\le i<j\le s}\lambda_i\centerdot\lambda_j .$ Similarly, we consider the $i$-th power $\Lambda_{\lambda}^i$ of $\Lambda_{\lambda}$ in the $s\times s$ matrix semiring $M_s(\mathbb{N}[\lambda_1,\cdots,\lambda_s])$ over $\mathbb{N}[\lambda_1,\cdots,\lambda_s]$ and then define $\|\Lambda_{\lambda}^i\| := \|\Lambda_{\lambda}^i\|_{_{\lambda}}= (\Lambda_{\lambda}^i)\varphi_{\lambda}$.
Note that $\dim_RS_n(\lambda,R)=g_1(\lambda)+2\|\Lambda_{\lambda}\|$ (see Lemma \ref{lem21}(1)).

The triangular divisor matrix $\Lambda_{\lambda}$ of $\lambda$ depends upon the ordering of $\lambda_1,\cdots,\lambda_s$. Nevertheless, we will show that the norms $\|-\|$ of triangular divisor matrices are independent of the choices of orderings.

\begin{Prop}\label{equ-matrix} Let $\lambda,\mu\in P(s,n)$. Then

$(1)$ $\|\Lambda^i_{\lambda}\|=g_{i+1}(\lambda)$ for $1\le i\le s-1$.

$(2)$ $\lambda\sim \mu$ if and only if $\|\Lambda^i_{\lambda}\| = \|\Lambda^i_{\mu}\|$
for all $1\le i\le s-1$, that is,
the norms of $\Lambda^i_{\lambda}$ and $ \Lambda^i_{\mu}$ are equal for all $1\le i\le s-1$.
\end{Prop}

{\it Proof.} (1) We check $\|\Lambda^i_{\lambda}\|=g_{i+1}(\lambda)$ for $1\le i\le s-1$.  Since the summands of
$\|\Lambda^i_{\lambda}\|= \|\Lambda^{i-1}_{\lambda}\Lambda_{\lambda}\|$ are of the form $\lambda_{j_1}\centerdot\cdots\centerdot\lambda_{j_i}\centerdot\lambda_{j_{i+1}}$ with $j_1<\cdots<j_i<j_{i+1}$, these are precisely the summands of $g_{i+1}(\lambda)$. Hence $\|\Lambda^i_{\lambda}\|= g_{i+1}(\lambda)$.

(2) Clearly, $g_1(\lambda)=g_1(\mu)=n$. Hence $g_j(\lambda)=g_{j}(\mu)$ for all $1\le j \le s$ if and only if $\|\Lambda^i_{\lambda}\| = \|\Lambda^i_{\mu}\|$
for all $1\le i\le s-1$. $\square$

\smallskip
An immediate consequence of Proposition \ref{equ-matrix} is that if $\{\mu_1, \cdots,\mu_s\}$ is a permutation of $\{\lambda_1,\cdots,\lambda_s\}$ then $\|\Lambda^i_{\lambda}\| = \|\Lambda^i_{\mu}\|$ for all $1\le i\le s-1$ because the $i$-th elementary symmetric polynomial $g_i(\lambda)$ does not depend on the ordering of its variables $\lambda_i$.

Another consequence is the following result which is quite useful for deciding whether two partitions are polynomial equivalent.

\begin{Koro}\label{equ-perm} Let $\lambda,\mu\in P(s,n)$ with the corresponding triangular divisor matrices $\Lambda_{\lambda}=(d_{ij})$ and $\Lambda_{\mu}=(c_{ij})$, respectively. If the multisets $\{d_{ij}\mid 1\le i<j\le s\}$ and $\{c_{ij}\mid 1\le i<j\le s\}$ are equal, then $\lambda\sim \mu$.
\end{Koro}

{\it Proof.} That the multisets $\{d_{ij}\mid 1\le i<j\le s\}$ and $\{c_{ij}\mid 1\le i<j\le s\}$ are equal means $g_2(\lambda)=g_2(\mu)$. Since any summand $ \lambda_{j_1}\centerdot\lambda_{j_2}\centerdot\cdots\centerdot\lambda_{j_i}$ of $g_i(\lambda)$ for $i\ge 3$ is a product $ d_{j_1,j_2}\centerdot d_{j_2,j_3}\centerdot\cdots\centerdot d_{j_{i-1},j_i}$ (in the canonical quarter-ring $(\mathbb{N},+,\centerdot)$ of $i-1$ elements of the multiset $\{d_{ij}\mid 1\le i<j\le s\}$ with the increasing indices, the summation of all such products of elements in the multiset of $\lambda$ equals the one in the multiset of $\mu$. Hence $g_i(\lambda)=g_i(\mu)$ for $2\le i\le s$. Clearly, $g_1(\lambda)=g_1(\mu)$. Thus $\lambda\sim \mu.$ $\square$

Note that the converse of Lemma \ref{equ-perm} fails. For example, $\lambda=(12,4,3,1)\sim \mu=(10,5,3,2)$, but the corresponding multisets for $\lambda$ and $\mu$ are $\{4,3,1,1,1,1\}$ and $\{5,1,2,1,1,1\}$, respectively, they are clearly different. This example also shows that if the common part $3$ is removed from $\lambda$ and $\mu$, then the resulting partitions $(12,4,1)$ and $(10,5,2)$ are no longer polynomial equivalent.

\begin{Koro} Let $\lambda=(\lambda_1,\cdots,\lambda_s), \mu=(\mu_1,\cdots,\mu_s) \in P(s,n)$. If $\lambda_i\centerdot\lambda_j=d=\mu_p\centerdot\mu_q$ for all $1\le i\ne j\le s$ and $1\le p\ne q\le s$, then $\lambda\sim \mu.$ In particular, if $\lambda_1,\cdots,\lambda_s$ are pairwise coprime and if $\mu_1,\cdots,\mu_s$ are pairwise coprime, then $\lambda\sim \mu$.
\label{equ-coprime}
\end{Koro}

{\it Proof.} Under the assumption, the triangular divisor matrices $\Lambda_{\lambda}$ and $\Lambda_{\mu}$ are equal. Hence Corollary \ref{equ-coprime} follows immediately from Corollary \ref{equ-perm}. $\square$

By the definition of partition polynomials $\epsilon_{\lambda}(x)$, we have the following useful fact, which implies that, when considering $\lambda\sim \mu$, we can always assume that $\lambda_1,\cdots,\lambda_s$ are coprime
and that $\mu_1,\cdots,\mu_s$ are coprime. By $d\lambda$ we means $(d\lambda_1, d\lambda_2,\cdots, d\lambda_s).$

\begin{Lem}\label{equ-mult} Let $\lambda, \mu\in P(s,n)$. Then $\lambda \sim \mu$ in $P(s,n)$ if and only if $d\lambda \sim d\mu$ in $P(s,dn)$ for some integer $d\ge 1$ if and only if $d\lambda \sim d\mu$ in $P(s,dn)$ for all integers $d\ge 1.$
\end{Lem}

{\it Proof.} The coefficients of the monic partition polynomial $\epsilon_{d\lambda}(x)$ are given by $g_i(d\lambda)/g_s(d\lambda)$. As we know, $g_i(d\lambda)=dg_i(\lambda)$ for all $i$. Thus $\epsilon_{d\lambda}(x)=\epsilon_{\lambda}(x)$. $\square$

\begin{Lem} Let $\lambda, \mu\in P(s,n)$ and $1\le m\in \mathbb{N}$ such that $\{\lambda_j\centerdot m \mid 1\le i\le s\} = \{\mu_j\centerdot m\mid 1\le j\le s\}$ as multisets.  Then $\lambda \sim \mu$ if and only if $(\lambda,m)\sim (\mu,m)$ in $P(s+1,n+m)$. In particular, suppose $\lambda_j\centerdot m = d = \mu_j\centerdot m$ for $1\le j\le s$. Then $\lambda \sim \mu$ if and only if $(\frac{\lambda}{d},\frac{m}{d})\sim (\frac{\mu}{d},\frac{m}{d})$ in $P(s+1, \frac{n+m}{d})$.
\label{equ-1}
\end{Lem}

{\it Proof.} The last conclusion in Lemma \ref{equ-1} follows from Lemma \ref{equ-mult}. We prove the other one. Let $\bar{\lambda}=(\lambda,m)\in P(s+1, n+m)$. Since the multisets $\{\lambda_j\centerdot m \mid 1\le i\le s\}$ and $\{\mu_j\centerdot m\mid 1\le j\le s\}$ are equal, the summation $\sum_{j_1<j_2 < \cdots < j_{i-1}}(\lambda_{j_1} \centerdot m) \centerdot\cdots\centerdot(\lambda_{j_{i-1}}\centerdot m)$ is equal to $\sum_{j_1<j_2<\cdots<j_{i-1}}(\mu_{j_1}\centerdot m) \centerdot\cdots\centerdot(\mu_{j_{i-1}}\centerdot m).$ In particular, $g_{s+1}(\bar{\lambda})=g_{s+1}(\bar{\mu})$. Clearly, for $1\le i\le s$, we have
$$g_i(\bar{\lambda})= g_i(\lambda) + \sum_{1\le j_i<\cdots<j_{i-1}\le s}\lambda_{j_1}\centerdot\cdots \centerdot\lambda_{j_{i-1}}\centerdot m = g_i(\lambda) + \sum_{1\le j_i<\cdots<j_{i-1}\le s}(\lambda_{j_1}\centerdot m)\cdots \centerdot (\lambda_{j_{i-1}}\centerdot m).$$ Thus, for $1\le i\le s$, $g_i(\bar{\lambda})= g_i(\bar{\mu})$ if and only if $g_i(\lambda)=g_i(\mu)$. This shows Lemma \ref{equ-1}. $\square$.

\medskip
A special case of Lemma \ref{equ-1} is that $\lambda\sim \mu$ if and only if $(\lambda,1)\sim (\mu,1)$. A slightly generalized form of the last statement of Lemma \ref{equ-1} is as follows.

\begin{Lem} Assume $\lambda \sim \mu$ in $P(s,n)$ and $\gamma \sim \delta$ in $P(t, m)$. If
$\lambda_i\centerdot\gamma_j= \mu_p\centerdot\delta_q$ for all $1\le i,p\le s, 1\le j,q\le t$, then $(\lambda,\gamma)\sim (\mu,\delta)$ in $P(s+t,n+m).$
\label{equ-m}
\end{Lem}

{\it Proof.} By assumption, we suppose $d:=\lambda_i\centerdot\gamma_j= \mu_p\centerdot \delta_q$. Now, let $\bar{\lambda}=(\lambda,\gamma)$ and $\bar{\mu}=(\mu,\delta)$. We calculate $g_i(\bar{\lambda})$ and $g_i(\bar{\mu})$ for $1\le i\le s+t.$ Clearly, $g_1(\lambda,\gamma)=g_1(\mu,\delta)$ and $g_{s+t}(\lambda,\gamma)=g_{s+t}(\mu,\delta)$. We may assume $s\le t$ and consider the following cases.

(1) $1< i\le s\le t$.
$$\begin{array}{rl} g_i(\bar{\lambda})& =\sum_{1\le j_1<\cdots<j_i\le s}\lambda_{j_1}\centerdot \cdots \centerdot\lambda_{j_i} + \sum_{1\le k_1<\cdots<k_i\le t}\gamma_{k_1}\centerdot \cdots\centerdot \gamma_{k_i} \\ & + \sum_{(p,q): p+q=i, 1\le p,q\le i}\sum_{1\le j_1<\cdots<j_p\le s, 1\le k_{1}< \cdots< k_q\le t}\lambda_{j_1}\centerdot\cdots\centerdot\lambda_{j_p}\centerdot \gamma_{k_1}\centerdot\cdots\centerdot\gamma_{k_q}\\
& = g_i(\lambda)+g_i(\gamma) + \sum_{(p,q): p+q=i, 1\le p\le s, 1\le q\le t} dC^p_sC^q_t\\
& = g_i(\mu)+g_i(\delta) + \sum_{(p,q): p+q=i, 1\le p\le s, 1\le q\le t} dC^p_sC^q_t\\
& = g_i(\mu)+g_i(\delta) + \sum_{(p,q): p+q=i, 1\le p,q\le i}\sum_{1\le j_1<\cdots<j_p\le s, 1\le k_{1}< \cdots< k_q\le t}\mu_{j_1}\centerdot\cdots\centerdot\mu_{j_p}\centerdot \delta_{k_1}\centerdot\cdots\centerdot\delta_{k_q}\\
& = g_i(\bar{\mu}).
\end{array}$$

(2) $s+1\le i< s+t$.
In this case, we first assume $i\le t$. Then
$$\begin{array}{rl} g_i(\bar{\lambda})& = \sum_{1\le k_1<\cdots<k_i\le t}\gamma_{k_1}\centerdot \cdots\centerdot \gamma_{k_i} \\ & + \sum_{(p,q): p+q=i, 1\le p,q\le i}\sum_{1\le j_1<\cdots<j_p\le s, 1\le k_{1}< \cdots< k_q\le t}\lambda_{j_1}\centerdot\cdots\centerdot\lambda_{j_p}\centerdot \gamma_{k_1}\centerdot\cdots\centerdot\gamma_{k_q}\\
& = g_i(\gamma) + \sum_{(p,q): p+q=i, 1\le p\le s, 1\le q\le t} dC^p_sC^q_t\\
& = g_i(\delta) +  \sum_{(p,q): p+q=i, 1\le p\le s, 1\le q\le t} dC^p_sC^q_t\\
& = g_i(\delta) + \sum_{(p,q): p+q=i, 1\le p,q\le i}\sum_{1\le j_1<\cdots<j_p\le s, 1\le k_{1}< \cdots< k_q\le t}\mu_{j_1}\centerdot\cdots\centerdot\mu_{j_p}\centerdot \delta_{k_1}\centerdot\cdots\centerdot\delta_{k_q}\\
& = g_i(\bar{\mu}).
\end{array}$$

Next, we assume $t<i$. Then
$$\begin{array}{rl} g_i(\bar{\lambda})& =  \sum_{(p,q): p+q=i, 1\le p,q\le i}\sum_{1\le j_1<\cdots<j_p\le s, 1\le k_{1}< \cdots< k_q\le t}\lambda_{j_1}\centerdot\cdots\centerdot\lambda_{j_p}\centerdot \gamma_{k_1}\centerdot\cdots\centerdot\gamma_{k_q}\\
& =  \sum_{(p,q): p+q=i, 1\le p\le s, 1\le q\le t} dC^p_sC^q_t\\
& = \sum_{(p,q): p+q=i, 1\le p,q\le i}\sum_{1\le j_1<\cdots<j_p\le s, 1\le k_{1}< \cdots< k_q\le t}\mu_{j_1}\centerdot\cdots\centerdot\mu_{j_p}\centerdot \delta_{k_1}\centerdot\cdots\centerdot\delta_{k_q}\\
& = g_i(\bar{\mu}).
\end{array}$$
Thus $g_i(\bar{\lambda})=g_i(\bar{\mu})$ for all $1\le i\le s+t$. $\square$

Finally, we point out that, for $\lambda\in P(s,n)$, $g_{i+1}(\lambda)$ can be calculated graphically by defining a valued quiver $Q_{\lambda}$ and writing out all paths of length $i$, and then summarising all evaluations of the paths in $(\mathbb{N}, +, \centerdot)$. For $i=0$, we understand $g_1(\lambda)=\sum_{i=1}^s\lambda_i=n$.
The quiver $Q_{\lambda}$ has the vertex set $\{1, 2, \cdots, s\}$, and for $i<j$, there is an arrow from $i$ to $j$ valued $d_{ij}$. For instance, $\lambda=(12,4,3,1)$, the quiver $Q_{\lambda}$ is
$$ 
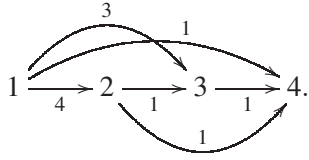
\xymatrix{
1\ar[r]_-{4}\ar@/^2pc/[rr]^{3}\ar@/^1.5pc/[rrr]^(0.6){1} & 2\ar[r]_{1}\ar@/_2pc/[rr]^{1} &3 \ar[r]_-{1} & 4.\\
}$$
Thus $g_4(\lambda)=4\centerdot 1 \centerdot 1=1, g_3(\lambda)=4\centerdot 1 + 4\centerdot 1 + 3\centerdot 1 + 1\centerdot 1= 4,$ $g_2(\lambda)=4+3+1+1+1+1=11, g_1(\lambda) =20$.

\subsection{Isomorphisms of invariant matrix algebras}
This section is devoted to a proof of Theorem \ref{iso-m}. We first show an auxiliary lemma.

\begin{Lem}\label{iso}
Let $R$ be a division ring of characteristic $p\ge 0$, and let $\sigma\in \Sigma_n$ be of cycle type $(\lambda_1,\lambda_2,\cdots,\lambda_s)$ with $\lambda_i\ge 1$ for all $1\le i\le s$. Then

$(1)$ $S_n(\sigma,R)$ is semisimple if and only if $p\nmid \lambda_i$ for $1\le i\le s$.

$(2)$ If $R$ is an algebraically closed field and $S_n(\sigma,R)$ is semisimple, then $$S_n(\sigma,R)\simeq R^{\oplus h_1(\sigma)}\times M_2(R)^{\oplus h_2(\sigma)}\times \cdots \times M_s(R)^{\oplus h_s(\sigma)}$$ as algebras, where $R^{\oplus h_1(\sigma)}$ denotes the direct product of $h_1(\sigma)$ copies of the ring $R$.
\end{Lem}

{\it Proof.} (1) Clearly, $M_n(R)$ is a simple artinian ring if $R$ is a division ring. Let $\rad(A)$ denote the Jacobson radical of a ring $A$. If $S_n(\sigma, R)$ is semisimple, that is $\rad\big(S_n(\sigma, R)\big)=0$, then $\rad\big(f_iS_n(\sigma, R)f_i\big)$ = $f_i\rad\big(S_n(\sigma, R)\big)f_i=0$ for all $1\le i\le s$. This means that all $f_iS_n(\sigma, R)f_i$ are semisimple. By Lemma \ref{lem2.4}(1), the group algebra $R[C_{\lambda_i}]$ is semisimple. Hence $p\nmid \lambda_i$ for all $i$. Conversely, if $p\nmid \lambda_i$ for all $i$, then $p$ does not divide the least common multiple of $\lambda_1, \lambda_2, \cdots, \lambda_s$, that is, the order of $G=\langle \sigma\rangle$ is invertible in $R$. In this case, $S_n(\sigma, R)$ is semisimple by \cite[Theorem 1.15, p.15]{Mont}, which says that if a finite group acts on a semisimple artinian ring with its order invertible in the ring, then the  ring of invariants is semisimple.

(2) Let $c_{\sigma}:=e_{1,(1)\sigma}+e_{2,(2)\sigma}+\cdots+e_{n,(n)\sigma}$ be the permutation matrix in $M_n(R)$ corresponding to $\sigma$.
Then $c_{\sigma}={\rm diag}\{c_{\sigma_1},c_{\sigma_2},\cdots,c_{\sigma_s}\}$, where $c_{\sigma_i}\in M_{\lambda_i}(R)$ is defined similarly.
For $1\le i\le s$, the eigenvalues of $c_{\sigma_i}$ are distinct and denoted by $V_i$.
Hence the Jordan canonical form of $c_{\sigma_i}$ has only Jordan blocks of order $1$, namely the diagonal matrix ${\rm diag}(1,\omega_i,\cdots,\omega_i^{\lambda_i-1})$, where $\omega_i$ is a primitive root of the unity of order $\lambda_i$ in $R$. This implies that
the Jordan canonical form of $c_{\sigma}$ is the diagonal matrix $$c:={\rm diag}(1,\omega_1,\cdots,\omega_1^{\lambda_1-1},\cdots,1,\omega_i,\cdots,\omega_i^{\lambda_i-1}, \cdots, 1,\omega_s,\cdots,\omega_s^{\lambda_s-1}).$$
For $1\le i\le s$, we write $H_i(\sigma) = \{r_{i1},r_{i2},\cdots,r_{ih_i(\sigma)}\}$ (see the notation above Definition \ref{character}).
Then the algebraic and geometric multiplicities of the eigenvalue $r_{ij}$ of $c_{\sigma}$ are just $i$ for $1\le j\le h_{i}(\sigma)$.
Thus $c$ is similar to
$$d:={\rm diag}(r_{11}I_1,r_{12}I_1,\cdots,r_{1\, h_1(\sigma)}I_1,
\cdots,r_{i1}I_i,r_{i2}I_i,\cdots,r_{i\, h_i(\sigma)}I_i,\cdots,
r_{s1}I_s,r_{s2}I_s,\cdots,r_{s\, h_s(\sigma)}I_s)$$
with $\sum_{i=1}^{s}i\,h_i(\sigma)=n.$ Recall that $I_i$ denotes the identity matrix in $M_i(R)$.
Hence the matrix $c_{\sigma}$ is similar to the matrix $d$.
By Lemma \ref{lem4.3} (see also \cite[Lemma 2.1]{Xi2021}), $S_n(\sigma,R)=S_n(c_{\sigma},R):=\{x\in M_n(R)\mid xc_{\sigma}=c_{\sigma}x\}\simeq S_n(d,R):=\{x\in M_n(R)\mid xd=dx\}$ as algebras.
Note that $r_{ij}\neq r_{kl}$ for $(i,j)\ne (k,l)$ with $1\le i,k\le s,1\le j\le h_{i}(\sigma),1\le l\le h_{k}(\sigma)$. It follows from  \cite[Lemma 2.2]{Xi2021} that the centralizer algebra $S_n(d,R)$ is isomorphic to $ R^{\oplus h_1(\sigma)}\times M_2(R)^{\oplus h_2(\sigma)}\times \cdots \times M_s(R)^{\oplus h_s(\sigma)}$.
Thus $S_n(\sigma,R)\simeq R^{\oplus h_1(\sigma)}\times M_2(R)^{\oplus h_2(\sigma)}\times \cdots \times M_s(R)^{\oplus h_s(\sigma)}$ as algebras.
$\square$

Lemma \ref{iso}(2) shows that the multiplicity of the matrix algebra $M_i(R)$ in a Wedderburn decomposition of the semisimple algebra $S_n(\lambda,R)$ is given by $h_i(\lambda)$.

Recall that two rings $R$ and $S$ are said to be \emph{Morita equivalent} if the two categories $R\Modcat$ and $S\Modcat$ (or equivalently $R\modcat$ and $S\modcat$) are equivalent.

\medskip
{\bf Proof of Theorem \ref{iso-m}.} (1) follows from Lemma \ref{iso}(1).

(2) Assume that $R$ is an algebraically closed field and that $S_n(\sigma,R)$ and $S_m(\tau,R)$ are semisimple. By Lemma \ref{iso}(2), we have $$S_n(\sigma,R)\simeq R^{\oplus h_1(\sigma)}\times M_2(R)^{\oplus h_2(\sigma)}\times \cdots \times M_s(R)^{\oplus h_s(\sigma)}
\mbox{ and }
S_m(\tau,R)\simeq R^{\oplus h_1(\tau)}\times M_2(R)^{\oplus h_2(\tau)}\times \cdots \times M_t(R)^{\oplus h_t(\tau)}.$$

(i) It follows from $h_i(\sigma)\ge 0,h_i(\tau)\ge 0$ for $1\le i\le s$, $h_s(\sigma)=\lambda_1\centerdot\cdots\centerdot\lambda_s\neq 0$ and $h_t(\tau)\neq 0$ that
$$R^{\oplus h_1(\sigma)}\times M_2(R)^{\oplus h_2(\sigma)}\times \cdots \times M_s(R)^{\oplus h_s(\sigma)}
\simeq
R^{\oplus h_1(\tau)}\times M_2(R)^{\oplus h_2(\tau)}\times \cdots \times M_t(R)^{\oplus h_t(\tau)}$$
as algebras if and only if $s=t$ and $h_i(\sigma)=h_i(\tau)$ for all $1\le i\le s.$

Suppose $S_n(\sigma,R)\simeq S_m(\tau,R)$, that is, we assume that $s=t$ and
$h_i(\sigma)=h_i(\tau)$ for all $i$. We show that $m=n$ and $\epsilon_{\lambda}(x)=\epsilon_{\mu}(x)$. In fact, since $\sum_{i=1}^{s}i\, h_i(\sigma)=n$ and $\sum_{i=1}^{t}i\, h_i(\tau)=m,$ we have $m=n$. By the definitions of $h_s(\sigma)$ and $h_s(\tau)$, we get $h_s(\sigma)=g_s(\sigma)$ and $h_s(\tau)=g_s(\tau)$. It follows from $h_s(\sigma)=h_s(\tau)$ that $g_s(\sigma)=g_s(\tau)$.
Now, by the definitions $h_{s-1}(\sigma)$ and $h_{s-1}(\tau)$ and Lemma \ref{in-ex}(1), we have $h_{s-1}(\sigma)=g_{s-1}(\sigma)-C_{s}^{s-1}g_{s}(\sigma)$ and $h_{s-1}(\tau)=g_{s-1}(\tau)-C_{s}^{s-1}g_{s}(\tau)$.
It then follows from $h_{s-1}(\sigma)=h_{s-1}(\tau)$ and $g_s(\sigma)=g_s(\tau)$ that $g_{s-1}(\sigma)=g_{s-1}(\tau)$.
Continuing this procedure, we get $g_{i}(\sigma)=g_{i}(\tau)$ for $1\le i\le s-2$.
Thus, if $h_i(\sigma)=h_i(\tau)$ for $1\le i\le s$, then $g_i(\sigma)=g_i(\tau)$ for $1\le i\le s$ and $\epsilon_{\lambda}(x)=\epsilon_{\mu}(x)$.

Conversely, if $m=n$ and $\epsilon_{\lambda}(x)=\epsilon_{\mu}(x)$, then $s=t$ and $g_s(\tau)g_i(\sigma)=g_s(\sigma)g_i(\tau)$ for $1\le i\le s-1$. As $g_1(\sigma)=n=m=g_1(\mu)$, we get $g_s(\sigma)=g_s(\tau)$. Thus it follows from $\epsilon_{\lambda}(x)=\epsilon_{\mu}(x)$ that $g_i(\sigma) =g_i(\tau)$ for all $1\le i\le s$. By Lemma \ref{in-ex}(2), we obtain $h_i(\sigma)=h_i(\tau)$ for all $i$.
Hence $S_n(\sigma,R)\simeq S_m(\tau,R)$.

(ii) Since $S_n(\sigma,R)\simeq R^{\oplus h_1(\sigma)}\times M_2(R)^{\oplus h_2(\sigma)}\times \cdots \times M_s(R)^{\oplus h_s(\sigma)}$ by Lemma \ref{iso}(2), the basic algebra $B(\sigma)$ of
$S_n(\sigma,R)$ is isomorphic to $R^{\oplus h_1(\sigma)}\times R^{\oplus h_2(\sigma)}\times\cdots\times R^{\oplus h_s(\sigma)}$. As is known, $S_n(\sigma,R)$ and $B(\sigma)$ are always Morita equivalent, while $S_n(\sigma,R)$ and $S_m(\tau,R)$ are Morita equivalent if and only if their basic algebras $B(\sigma)$ and $B(\tau)$ are isomorphic. Clearly, $B(\sigma)\simeq B(\tau)$ if and only if
$\sum_{i=1}^{s}h_i(\sigma)=\sum_{i=1}^{t}h_i(\tau)$.
By Lemma \ref{in-ex}(2), $\sum_{i=1}^{s}h_i(\sigma)=\sum_{i=1}^{t}h_i(\tau)$
 if and only if the matrices $c_{\sigma}$ and $c_{\tau}$ have the same number of distinct eigenvalues in $R$ if and only if
$\sum_{i=1}^{s}(-1)^{i+1}g_i(\sigma)=\sum_{i=1}^{t}(-1)^{i+1}g_i(\tau)$.
Therefore $S_n(\sigma,R)$ and $S_m(\tau,R)$ are Morita equivalent if and only if $\sum_{i=1}^s(-1)^{i+1}g_i(\sigma)=\sum_{i=1}^t(-1)^{i+1} g_i(\tau)$ if and only if $(-1)^{s-1}g_s(\sigma)\epsilon_{\lambda}(1)=(-1)^{t-1}g_t(\tau)\epsilon_{\mu}(1)$. Note that $g_s(\sigma)$ is the greatest common divisor of $\lambda_1,\cdots,\lambda_s$. $\square$

\medskip
Theorem \ref{iso-m} tells us that the foregoing consideration on the polynomial equivalence $\lambda\sim \mu$ in $P(s,n)$ provides us with many isomorphisms of invariant matrix algebras. For example, $S_{38}\big((17, 11, 8, 2), \mathbb{C}\big) \simeq  S_{38}\big((17, 11, 6, 4), \mathbb{C}\big)$ by repeatedly applying Corollary \ref{equ-1} to the polynomial equivalence $(8,2) \sim (6,4)$ in $P(2,10)$. Also, for any prime number $p\ge 5$ and two-part partitions $\lambda$ and $\mu$ of $p$, it follows from Theorem \ref{iso-m} that $S_p(\lambda,\mathbb{C})\simeq S_p(\mu,\mathbb{C})$. For instance, $S_5((4,1),\mathbb{C})\simeq S_5((3,2),\mathbb{C})\simeq \mathbb{C}\times\mathbb{C}\times\mathbb{C}\times M_2(\mathbb{C})$.

Remark that the condition $m=n$ in Theorem \ref{iso-m}(2)(i) is needed. For example, both $\lambda=(4,2,2)$ and $\mu=(2,1,1)$ have the same partition polynomial $x^2-3x+4$, but $S_8(\lambda,R)\not\simeq S_4(\mu,R)$ for any ring $R$.

\medskip
As a consequence of Theorem \ref{iso-m} and Corollary \ref{equ-coprime}. we have the corollary immediately.
\begin{Koro}
Let $R$ be an algebraically closed field of characteristic $p\ge 0$, and let $\sigma\in \Sigma_n$ and $\tau \in \Sigma_m$ be of cycle types $(\lambda_1,\lambda_2,\cdots,\lambda_s)$ and $(\mu_1,\mu_2,\cdots,\mu_t)$, respectively. Assume that $p\nmid \lambda_i$ and $p\nmid\mu_j$ for $1\le i\le s$ and $1\le i\le t$. If
$\lambda_1,\cdots,\lambda_s$ are pairwise coprime and  if $\mu_1,\cdots,\mu_t$ are pairwise coprime, then

$(1)$ $S_n(\sigma,R)\simeq S_m(\tau,R)$  if and only if $n=m$ and $s=t$.

$(2)$ $S_n(\sigma,R)$ and $S_m(\tau,R)$ are Morita equivalent if and only if $n-s=m-t$.
\end{Koro}

\medskip
Finally, we mention that there are a few issues worthy of further consideration, which have
directly been suggested in this section but not answered.
We formalize some of them as the following questions.

\medskip
{\bf Questions.} (1) When are $S_n(\lambda, R)$ and $S_n(\mu,R)$ isomorphic (or derived equivalent) for $\lambda,\mu\in P(s,n)$ and $R$ a unitary ring?

(2) Given a polynomial $f(x)=x^{s-1} -a_{s-2}x^{s-2}+\cdots +(-1)^{s-i-1}a_{i}x^i+\cdots + (-1)^{s-2}a_{1}x +(-1)^{s-1}a_0 \in \mathbb{Z}[x]$ with $s\ge 3$ and $a_i\ge C_s^{i+1}$ for $1\le i\le s-2$, is there a partition $\lambda\in P(s,n)$ for some $n$ such that $\epsilon_{\lambda}(x)=f(x)$?  In other words, what are the necessary and sufficient conditions for a monic polynomial $f(x)$ to be a partition polynomial?

(3) Are $\epsilon_{\lambda}(x)\epsilon_{\mu}(x)$ and $\epsilon_{\lambda}(\epsilon_{\mu}(x))$ again partition polynomials? Or more generally, how do the operations of partitions correspond to the ones of partition polynomials?

(4) Let $I(s,n)$ be the set of equivalence classes of $P(s,n)$ with respect to the polynomial equivalence relation. Then $I(s,n)$ is the union of distinct equivalence classes. Let $E(s,n,j)$ be the set of those equivalence classes that contain $j$ elements of $P(s,n)$. Define $p(s,n)$, $i(s,n)$ and $e(s,n,j)$ to be the numbers of elements in $P(s,n)$, $I(s,n)$ and $E(s,n,j)$, respectively. What are the generating functions of $i(s,n)$ and $e(s,n,j)$?

Note that the generating function for $p(s,n)$ is
$ \sum_{n=0}^{\infty}p(s,n)q^n = \frac{q^s}{(1-q)(1-q^2)\cdots(1-q^s)}$ (see \cite[Chapter 6]{ae}).

(5) Describe partitions $\lambda\in P(s,n)$ that are self-equivalent, that is, if $\mu\in P(s,n)$ such that $\mu \sim \lambda$, then $\mu=\lambda$. For instance, $\lambda=(4,4,1)$ is self-equivalent with $\epsilon_{(4,4,1)}(x)=(x-3)^2$, but $\mu=(5,2,2)$ is not self-equivalent.

(6) For which partitions $\lambda$ are $\epsilon_{\lambda}(x)$ irreducible polynomials over $\mathbb{Q}$?

(7) What are the necessary and sufficient conditions for $\lambda\sim \mu$ in $P(s,n)$ in terms of the conjugate partitions of $\lambda$ and $\mu$?

\medskip
{\bf Acknowledgement.}
The research work of both authors was partially supported by the National Natural Science Foundation of China (12031014). The author CCX thanks Professors Hourong Qin from Nanjing University for discussions on greatest common divisor matrices, and Yong Shao from Northwest University in Xi'an for conversations on semirings. Also, the authors are very grateful to the referee for very detailed comments, helpful suggestions and interesting questions which improve both the presentation and some results of the paper.

\medskip
{\footnotesize

Changchang Xi, School of Mathematical Sciences, Capital Normal University, 100048 Beijing, P.R. China

{\tt Email: xicc@cnu.edu.cn}

\medskip
Jinbi Zhang, School of Mathematical Sciences, Peking University, 100871 Beijing, P.R. China

{\tt Email: zhangjb@cnu.edu.cn}
}
\end{document}